\documentclass[aos,MSNbibl,nameyear,dvips]{arximspdf}
\usepackage{dcolumn}

\doi{10.1214/12-AOS990} 
\volume{40}
\issue{3}
\pubyear{2012}
\firstpage{1816}
\lastpage{1845}

\makeatletter
\newcolumntype{d}[1]{D{.}{.}{#1}}
\newtheorem{theorem}{Theorem}
\makeatother

\begin{document}
\begin{frontmatter}

\title{Semiparametric theory for causal mediation analysis:
Efficiency bounds, multiple robustness and~sensitivity~analysis}
\runtitle{Estimation of the mediation functional}

\begin{aug}
\author[A]{\fnms{Eric J.} \snm{Tchetgen Tchetgen}\corref{}\thanksref{t1}\ead[label=e1]{etchetge@hsph.harvard.edu}}
\and
\author[B]{\fnms{Ilya} \snm{Shpitser}\ead[label=e2]{ishpitse@hsph.harvard.edu}}
\runauthor{E. J. Tchetgen Tchetgen and I. Shpitser}
\thankstext{t1}{Supported by NIH Grant R21ES019712.}
\affiliation{Harvard School of Public Health}
\address[A]{Departments of Epidemiology and Biostatistics\\
Harvard School of Public Health\\
677 Huntington Avenue\\
Boston, Massachusetts 02115\\
USA\\
\printead{e1}} 
\address[B]{Department of Epidemiology\\
Harvard School of Public Health\\
677 Huntington Avenue\\
Boston, Massachusetts 02115\\
USA\\
\printead{e2}}
\end{aug}

\received{\smonth{3} \syear{2011}}
\revised{\smonth{3} \syear{2012}}

%
\begin{abstract}
While estimation of the marginal (total) causal effect of a point exposure
on an outcome is arguably the most common objective of experimental and
observational studies in the health and social sciences, in recent years,
investigators have also become increasingly interested in mediation
analysis. Specifically, upon evaluating the total effect of the exposure,
investigators routinely wish to make inferences about the direct or indirect
pathways of the effect of the exposure, through a mediator
variable or not, that occurs subsequently to the exposure and prior to
the outcome.
Although powerful semiparametric methodologies have been developed to
analyze observational studies that produce double robust and highly
efficient estimates of the marginal total causal effect, similar
methods for
mediation analysis are currently lacking. Thus, this paper develops
a~general semiparametric framework for obtaining inferences about so-called
marginal natural direct and indirect causal effects, while appropriately
accounting for a large number of pre-exposure confounding factors for the
exposure and the mediator variables. Our analytic framework is particularly
appealing, because it gives new insights on issues of efficiency and
robustness in the context of mediation analysis. In particular, we propose
new multiply robust locally efficient estimators of the marginal natural
indirect and direct causal effects, and develop a~novel double robust
sensitivity analysis framework for the assumption of ignorability of the
mediator variable.
\end{abstract}

%

\end{frontmatter}

\section{Introduction}

The evaluation of the total causal effect of a given point exposure,
treatment or intervention on an outcome of interest is arguably the most
common objective of experimental and observational studies in the
fields of
epidemiology,\vadjust{\goodbreak} biostatistics and in the social sciences. However, in recent
years, investigators in these various fields have become increasingly
interested in making inferences about the direct or indirect pathways
of the
exposure effect, through a mediator variable or not, that occurs
subsequently to the exposure and prior to the outcome. Recently, the
counterfactual language of causal inference has proven particularly useful
for formalizing mediation analysis. Indeed, causal inference offers a formal
mathematical framework for defining varieties of direct and indirect
effects, and for establishing necessary and sufficient identifying
conditions of these effects. A notable contribution of causal inference to
the literature on mediation analysis is the key distinction drawn between
so-called controlled direct effects versus natural direct effects. In words, the controlled direct effect refers to the
exposure effect that arises upon intervening to set the mediator to a fixed
level that may differ from its actual observed value [\citet{robins1}, \citet{Pearl}, \citet{Rob2003}]. In contrast, the natural
(also known as
pure) direct effect captures the effect of the exposure when one intervenes
to set the mediator to the (random) level it would have been in the absence
of exposure [\citet{robins1}, \citet{Pearl}]. As noted by
\citet{Pearl}, controlled direct and indirect effects are particularly
relevant for policy making, whereas natural direct and indirect effects are
more useful for understanding the underlying mechanism by which the exposure
operates.
In fact,
natural direct and indirect effects combine to produce the exposure total
effect.

To formally define natural direct and indirect effects first requires
defining counterfactuals. We assume that for each level of a binary exposure
$E$, and of a mediator variable $M$, there exist a counterfactual
variable $Y_{e,m}$ corresponding to the outcome $Y$ had possibly contrary to fact the
exposure and mediator variables taken the value $(e,m)$. Similarly, for
$E=e$, we assume there exists a counterfactual variable $M_{e}$
corresponding to
the mediator variable had possibly contrary to fact the exposure variable
taken the value $e$. The current paper concerns the decomposition of the
total effect of $E$ on $Y$, in terms of natural direct and natural indirect
effects, which, expressed on the mean difference scale, is given by%
\begin{eqnarray}\label{total_effect}
\overbrace{\mathbb{E} (Y_{e=1}-Y_{e=0} )}^{\mathrm{total\ effect}}
&=&\mathbb{E} (
Y_{e=1,M_{e=1}}-Y_{e=0,M_{e=0}} ) \nonumber\hspace*{-25pt}\\[-8pt]\\[-8pt]
\qquad&=&\overbrace{\mathbb{E} (
Y_{e=1,M_{e=1}}-Y_{e=1,M_{e=0}} ) }^{\mathrm{natural\ indirect\ effect}}+
\overbrace{\mathbb{E} ( Y_{e=1,M_{e=0}}-Y_{e=0,M_{e=0}} ) }%
^{\mathrm{natural\ direct\ effect}}, \nonumber\hspace*{-25pt}
\end{eqnarray}
where $\mathbb{E}$ stands for expectation.

In an effort to account for confounding bias when estimating causal effects,
such as the average total effect (\ref{total_effect}) from
nonexperimental data, investigators routinely collect and adjust for in
data analysis, a large number of confounding factors. Because of the curse
of dimensionality, nonparametric methods\vadjust{\goodbreak} of estimation are typically not
practical in such settings, and one usually resorts to one of two
dimension-reduction strategies; either one relies on a model for the outcome
given exposure and counfounders, or alternately one relies on a model for
the exposure, that is, the propensity score. Recently, powerful semiparametric
methods have been developed to analyze observational studies that produce
so-called double robust and highly efficient estimates of the exposure total
causal effect [Robins (\citeyear{Robins3}), \citet{robins55}, \citet
{Robins4}, \citet{Tsiatis}] and similar methods have also been
developed to estimate controlled direct effects
[\citet{stijn}]. An important advantage of a double
robust method is that it carefully combines both of the aforementioned
dimension reduction strategies for confounding adjustment, to produce an
estimator of the causal effect that remains consistent and asymptotically
normal, provided at least one of the two strategies is correct, without
necessarily knowing which strategy is indeed correct [\citet
{vanderLaan}]. Unfortunately, similar methods for making semiparametric
inferences about marginal natural direct and indirect effects are currently
lacking. Thus, this paper develops a general semiparametric framework for
obtaining inferences about marginal natural direct and indirect effects on
the mean of an outcome, while appropriately accounting for a large
number of
confounding factors for the exposure and the mediator
variables.

Our semiparametric framework is particularly appealing, as it gives new
insight on issues of efficiency and robustness in the context of mediation
analysis. Specifically, in Section~\ref{sec2}, we adopt the sequential
ignorability
assumption of Imai, Keele and Tingley
(\citeyear{Imai2}) under which, in conjunction with the
standard consistency and positivity assumptions, we derive the efficient
influence function and thus obtain the semiparametric efficiency bound for
the natural direct and natural indirect marginal mean causal effects,
in the
nonparametric model $\mathcal{M}_{\mathrm{nonpar}}$ in which the
observed data
likelihood is left unrestricted. We further show that in order to conduct
mediation inferences in $\mathcal{M}_{\mathrm{nonpar}}$, one must
estimate at
least a subset of the following quantities:
\begin{longlist}[(iii)]
\item[(i)] the conditional expectation of the outcome given the mediator,
exposure and confounding factors;

\item[(ii)] the density of the mediator given the exposure and the
confounders;

\item[(iii)] the density of the exposure given the confounders.
\end{longlist}

Ideally, to minimize the possibility of modeling bias, one may wish to
estimate each of these quantities nonparametrically; however, as previously
argued, when as we assume throughout, we wish to account for numerous
confounders, such nonparametric estimates will likely perform poorly in
finite samples. Thus, in Section~\ref{sec2.3} we develop an alternative multiply
robust strategy. To do so, we propose to model (i), (ii) and (iii)
parametrically (or semiparametrically), but rather than obtaining mediation
inferences that rely on the correct specification of a specific subset of
these models, instead we carefully combine these three models to produce
estimators of\vadjust{\goodbreak} the marginal mean direct and indirect effects that remain
consistent and asymptotically normal (CAN) in a union model, where at least
one but not necessarily all of the following conditions hold:
\begin{longlist}[(b)]
\item[(a)] the parametric or semi-parametric models for the conditional
expectation of the outcome (i) and for the conditional density of the
mediator (ii) are correctly specified;

\item[(b)] the parametric or semiparametric models for the conditional
expectation of the outcome (i) and for the conditional density of the
exposure (iii) are correctly specified;

\item[(c)] the parametric or semiparametric models for the conditional
densities of the exposure and the mediator (ii) and (iii) are correctly
specified.
\end{longlist}

Accordingly, we define submodels $\mathcal{M}_{a}$, $\mathcal{M}_{b}$
and $%
\mathcal{M}_{c}$ of $\mathcal{M}_{\mathrm{nonpar}}$ corresponding to
models (a), (b) and (c) respectively. Thus, the proposed
approach is triply robust as it produces valid inferences about natural
direct and indirect effects in the union model $\mathcal{M}_{\mathrm{union}}=%
\mathcal{M}_{a}\cup\mathcal{M}_{b}\cup\mathcal{M}_{c}$. Furthermore, as
we later show in Section~\ref{sec2.3}, the proposed estimators are also locally
semiparametric efficient in the sense that they achieve the respective
efficiency bounds for estimating the natural direct and indirect
effects in $%
\mathcal{M}_{\mathrm{union}}$, at the intersection submodel $\mathcal{M}
_{a}\cap\mathcal{M}_{b} \cap \mathcal{M}_{c}=\mathcal{M}_{a}\cap
\mathcal{M}_{c}=\mathcal{M}_{a}\cap\mathcal{M}_{b}=\mathcal{M}_{b}\cap
\mathcal{M}_{c}\subset\mathcal{M}_{\mathrm{union}}\subset\mathcal
{M}_{\mathrm{%
nonpar}}$.

Section~\ref{sec3} summarizes a simulation study illustrating the
finite sample
performance of the various estimators described in Section~\ref{sec2},
and Section~\ref{sec4}
gives a real data application of these methods. Section~\ref{sec5}
describes a
strategy to improve the stability of the proposed multiply robust estimator
which directly depends on inverse exposure and mediator density weights,
when such weights are highly variable, and Section~\ref{sec6}
demonstrates the
favorable performance of two modified multiply robust estimators in the
context of such highly variable weights. In Section~\ref{sec7}, we
compare the
proposed methodology to the prevailing estimators in the literature. Based
on this comparison, we conclude that the new approach should generally be
preferred because an inference under the proposed method is guaranteed to
remain valid under many more data generating laws than an inference
based on
each of the other existing approaches. In particular, as we argue below the
approach of \citet{vanderlaan2} is not entirely satisfactory
because, despite producing a CAN estimator of the marginal direct effect
under the union model $\mathcal{M}_{a}\cup\mathcal{M}_{c}$ (and therefore
an estimator that is double robust), their estimator requires a correct
model for the density of the mediator. Thus, unlike the direct effect
estimator developed in this paper, the van der Laan estimator fails to be
consistent under the submodel $\mathcal{M}_{b}\subset\mathcal
{M}_{\mathrm{%
union}}$. Nonetheless, the estimator of van der Laan is in fact locally
efficient in model $\mathcal{M}_{a}\cup\mathcal{M}_{c}$, provided the model
for the mediator's conditional density is either known, or can be
efficiently estimated. This property is confirmed in a supplementary
online Appendix\vadjust{\goodbreak} [\citet{suppl}], where we also provide a general
map that relates the
efficient influence function for model $\mathcal{M}_{\mathrm{union}}$
to the
corresponding efficient influence function for model $\mathcal{M}_{a}\cup\mathcal
{M}%
_{c}$, assuming an arbitrary parametric or semiparametric model for the
mediator conditional density is correctly specified. In Section~\ref
{sec8}, we
describe a novel double robust sensitivity analysis framework to assess the
impact on inferences about the natural direct effect, of a departure from
the ignorability assumption of the mediator variable. We conclude with a
brief discussion.\vspace*{-2pt}

\section{The nonparametric mediation functional}\label{sec2}\vspace*{-2pt}

\subsection{Identification}\label{2.1}

Suppose i.i.d. data on $O=(Y,E,M,X)$ is collected for $n$ subjects.
Recall that $Y$ is an outcome of interest, $E$ is a binary exposure
variable, $M$ is a mediator variable with support $\mathcal{S}$, known to
occur subsequently to $E$ and prior to $Y$ and $X$ is a vector of
pre-exposure variables with support $\mathcal{X}$ that confound the
association between $(E,M)$ and $Y$. The overarching goal of this paper is
to provide some theory of inference about the fundamental functional of
mediation analysis which Judea Pearl calls ``the mediation causal formula''
[Pearl (\citeyear{Pearl2})] and which, expressed on the mean scale, is %
%
\begin{eqnarray}\label{main_functional}
\theta_{0}&=&\mathop{\int\hspace*{-4pt}\int}_{\mathcal{S\times X}}
\mathbb{E} ( Y|E=1,M=m,X=x )\nonumber\\[-8pt]\\[-8pt]
&&\hspace*{19pt}{}\times f_{M|E,X} ( m|E=0,X=x )
f_{X}(x)\,d\mu(m,x),\nonumber
\end{eqnarray}
$f_{M|E,X}$ and $f_{X}$ are respectively the conditional density of the
mediator $M$ given $(E,X)$ and the density of $X$, and $\mu$ is a
dominating measure for the distribution of $(M,X)$. Hereafter, to keep with
standard statistical parlance, we shall simply refer to $\theta_{0}$
as the
``mediation functional'' or ``M-functional'' since it is formally a functional
on the nonparametric statistical model $\mathcal{M}_{\mathrm{nonpar}}=\{
F_{O} (
\cdot ) \dvtx F_{O}\mbox{ unrestricted}\}$ of all regular laws $F_{O}$
of the
observed data $O$ that satisfy the positivity assumption given below;
that is, $\theta_{0}=\theta_{0} ( F_{O} ) \dvtx\mathcal{M}%
_{\mathrm{nonpar}}\rightarrow\mathcal{R}$, with $\mathcal{R}$ the real
line. The
functional $\theta_{0}$ is of keen interest here because it arises in the
estimation of natural direct and indirect effects as we describe next.
To do so, we make the consistency assumption.

\textit{Consistency:}
\begin{eqnarray*}
&&\mbox{if }E =e ,\mbox{ then }M_{e}=M\mbox{ w.p.1} \quad\mbox{and}\\
&&\mbox{if }E =e \mbox{ and }M=m, \mbox{ then }Y_{e,m}=Y\mbox{ w.p.1.}
\end{eqnarray*}

In addition, we adopt the sequential ignorability assumption of Imai, Keele and Tingley
(\citeyear{Imai2}) which states that for $e,e^{\prime}\in \{ 0,1 \} $.

\textit{Sequential ignorability:}%
\begin{eqnarray*}
&&\{ Y_{e^{\prime},m},M_{e} \} \perp\!\!\!\!\perp E|X, \label{A.2} \\
&&Y_{e^{\prime}m}\perp\!\!\!\!\perp M|E=e,X,\vadjust{\goodbreak}
\end{eqnarray*}
where $A\perp\!\!\!\!\perp B|C$ states that $A$ is independent of $B$
given $%
C$; paired with the following:

\textit{Positivity:}
\begin{eqnarray*}
f_{M|E,X} ( m|E,X ) & >&0\qquad\mbox{w.p.1 for each }m\in\mathcal{S}%
\quad\mbox{and} \\
f_{E|X} ( e|X ) & >&0\qquad\mbox{w.p.1 for each }e\in \{
0,1 \} .
\end{eqnarray*}

Then, under the consistency, sequential ignorability and positivity
assumptions, Imai, Keele and Tingley
(\citeyear{Imai2}) showed that
\begin{eqnarray}\label{ATE functional}
\theta_{0}& =&\mathbb{E} ( Y_{1,M_{0}} ) \quad\mbox{and } \nonumber\\
\delta_{e}& \equiv&\mathop{ \int}_{\mathcal{X}}\mathbb{E} ( Y|E=e,X=x ) f_{X}(x)\,d\mu(x) \nonumber\\
& =&\mathop{\int\hspace*{-4pt}\int}_{\mathcal{S\times X}}\mathbb{E}%
( Y|E=e,M=m,X=x )\\
&&\hspace*{17pt}{}\times f_{M|E,X} ( m|E=e,X=x ) f_{X}(x)\,d\mu
(m,x)  \nonumber\\
& =&\mathbb{E} ( Y_{e} ) =\mathbb{E} ( Y_{e,M_{e}} ) ,\qquad e=0,1, \nonumber
\end{eqnarray}
so that $\mathbb{E} ( Y_{1,M_{0}} ) $ and $\mathbb{E} (
Y_{e} ) $, $e=0,1$, are identified from the observed data, and so is
the mean natural direct effect $\mathbb{E} ( Y_{1,M_{0}} ) -\mathbb{%
E} ( Y_{0} ) =\theta_{0}-\delta_{0}$ and the mean natural
indirect effect $\mathbb{E} ( Y_{1} ) -\mathbb{E} (
Y_{1,M_{0}} ) =\delta_{1}-\theta_{0}$. For binary $Y$, one might
alternatively consider the natural direct effect on the risk ratio
scale $%
\mathbb{E} ( Y_{1,M_{0}} ) /\mathbb{E} ( Y_{0} ) =\theta
_{0}/\delta_{0}$ or on the odds ratio scale $ \{ \mathbb{E} (
Y_{1,M_{0}} ) \mathbb{E} ( 1-Y_{0} ) \} / \{
\mathbb{E} ( 1-Y_{1,M_{0}} ) \mathbb{E} ( Y_{0} ) \}
= \{ \theta_{0} ( 1-\delta_{0} ) \} / \{ \delta
_{0} ( 1-\theta_{0} ) \} $ and similarly defined natural
indirect effects on the risk ratio and odds ratio scales. It is
instructive to contrast the expression (\ref{main_functional})
for $\mathbb{E} ( Y_{1,M_{0}} )$ with the expression  (\ref{ATE functional})  for $e=1$ corresponding to $\mathbb{E} (
Y_{1} ) $, and to note that the two expressions bare a striking
resemblance except the density of the mediator in the first expression
conditions on the unexposed (with $E=0$), whereas in the second expression,
the mediator density is conditional on the exposed (with \mbox{$E=1$}). As we
demonstrate below, this subtle difference has remarkable implications for
inference.

\citet{Pearl} was the first to derive the M-functional $\theta
_{0}=\mathbb{E}%
( Y_{1,M_{0}} ) $ under a different set of assumptions. Others
have since contributed alternative sets of identifying assumptions. In this
paper, we have chosen to work under the sequential ignorability assumption
of \citet{Imai1}, \citet{Imai2}, but note that alternative related assumptions exist
in the literature [\citet{robins1}, \citet{Pearl}, \citet{vanderlaan2}, Hafeman and Vanderweele (\citeyear{Hafeman})]; however, we note that
\citet{Robins5} disagree with the label ``sequential
ignorability'' because its terminology has previously carried a~different
interpretation\vadjust{\goodbreak} in the literature. Nonetheless, the assumption entails two
ignorability-like assumptions that are made sequentially. First, given the
observed pre-exposure confounders, the exposure assignment is assumed
to be
ignorable, that is, statistically independent of potential outcomes and
potential mediators. The second part of the assumption states that the
mediator is ignorable given the observed exposure and pre-exposure
confounders. Specifically, the second part of the sequential ignorability
assumption is conditional on the observed value of the ignorable
treatment and the observed pretreatment confounders. We note that the second
part of the sequential ignorability assumption is particularly strong and
must be made with care. This is partly because it is always possible that
there might be unobserved variables that confound the relationship between
the outcome and the mediator variables, even upon conditioning on the
observed exposure and covariates. Furthermore, the confounders $X$ must all
be pre-exposure variables; that is, they must precede $E$. In fact,
\citet{Avin} proved that without additional assumptions, one cannot identify
natural direct and indirect effects if there are confounding variables that
are affected by the exposure, even if such variables are observed by the
investigator [also see Tchetgen Tchetgen and VanderWeele (\citeyear{TchVand2012})]. This implies that, similarly to the ignorability of the exposure
in observational studies, ignorability of the mediator cannot be established
with certainty, even after collecting as many pre-exposure confounders as
possible. Furthermore, as \citet{Robins5} point out, whereas
the first part of the sequential ignorability assumption could, in principle,
be enforced in a~randomized study, by randomizing $E$ within levels of $X;$
the second part of the sequential ignorability assumption cannot similarly
be enforced experimentally, even by randomization. And thus, for this latter
assumption to hold, one must entirely rely on expert knowledge about the
mechanism under study. For this reason, it will be crucial in practice to
supplement mediation analyses with a sensitivity analysis that accurately
quantifies the degree to which results are robust to a potential violation
of the sequential ignorability assumption. Later in the paper, we
develop a variety of sensitivity analysis techniques that allow the analyst to quantify the
degree to which his or her mediation analysis results are robust to a
potential violation of the sequential ignorability assumption.

\subsection{\texorpdfstring{Semiparametric efficiency bounds for $\mathcal{M}_{\mathrm{nonpar}}$}
{Semiparametric efficiency bounds for M nonpar}}\label{sec2.2}

In this section, we derive the efficient influence function for the
M-functional $\theta_{0}$ in $\mathcal{M}_{\mathrm{nonpar}}$. This
result is then
combined with the efficient influence function for the functional
$\delta_{e}$ [\citet{Robins6}, \citet{Hahn}] to obtain the efficient
influence function for the natural direct and indirect effects on the mean
difference scale. Thus, in the following, we shall use the efficient
influence function $S_{\delta e}^{\mathrm{eff},\mathrm{nonpar}} ( \delta_{e} ) $ of $%
\delta_{e}$ which is well known to be
\[
\frac{I(E=e)}{f_{E|X}(e|X)} \{ Y-\eta ( e,e,X ) \} +\eta
( e,e,X ) +\delta_{e},
\]
where for $e,e^{\ast}\in \{ 0,1 \} $, we define
\[
\eta ( e,e^{\ast},X ) =\int_{\mathcal{S}}\mathbb{E} (
Y|X,M=m,E=e ) f_{M|E,X} ( m|E=e^{\ast},X ) \,d\mu(m),
\]
so that $\eta ( e,e,X ) =\mathbb{E} ( Y|X,E=e ) $, $%
e=0,1$.

The following theorem is proved in the \hyperref[app]{Appendix}.

\begin{theorem}\label{teo1}
Under the consistency, sequential ignorability and
positivity assumptions, the efficient influence function of the M-functional
$\theta_{0}$ in model $\mathcal{M}_{\mathrm{nonpar}}$ is given by
\begin{eqnarray*}
&&S_{\theta_{0}}^{\mathrm{eff},\mathrm{nonpar}} ( \theta_{0} )\\
&&\qquad =s_{\theta_{0}}^{\mathrm{eff},\mathrm{nonpar}} ( O;\theta_{0})\\
 &&\qquad=\frac{I \{E=1 \} f_{M|E,X} ( M|E=0,X ) }{f_{E|X}(1|X)f_{M|E,X} (
M|E=1,X ) } \{ Y-\mathbb{E} ( Y|X,M,E=1 )  \} \\
&&\qquad\quad{}+\frac{I(E=0)}{f_{E|X}(0|X)} \{ \mathbb{E} ( Y|X,M,E=1 )
-\eta ( 1,0,X ) \} +\eta ( 1,0,X ) -\theta_{0},
\end{eqnarray*}
and the efficient influence function of the natural direct
and indirect effects on the mean difference scale in model $\mathcal{M}
_{\mathrm{nonpar}}$ are respectively given by
\begin{eqnarray*}
&&S_{\mathrm{NDE}}^{\mathrm{eff},\mathrm{nonpar}} (
\theta_{0},\delta_{0} )\\
&&\qquad = s_{\mathrm{NDE}}^{\mathrm{eff},\mathrm{nonpar}} ( O;\theta_{0},\delta_{0} )\\
&&\qquad =
S_{\theta_{0}}^{\mathrm{eff},\mathrm{nonpar}} ( \theta_{0} ) -S_{\delta
_{0}}^{\mathrm{eff},\mathrm{nonpar}} ( \delta_{0} ) \\
&&\qquad=\frac{I \{ E=1 \} f_{M|E,X} ( M|E=0,X ) }{%
f_{E|X}(1|X)f_{M|E,X} ( M|E=1,X ) } \{ Y-\mathbb{E} (
Y|X,M,E=1 )  \} \\
&&\quad\qquad{}+\frac{I(E=0)}{f_{E|X}(0|X)} \{ \mathbb{E} ( Y|X,M,E=1 )
-Y-\eta ( 1,0,X ) +\eta ( 0,0,X ) \} \\
&&\quad\qquad{}+\eta ( 1,0,X ) -\eta ( 0,0,X ) -\theta_{0}+\delta
_{0},
\end{eqnarray*}
  and
\begin{eqnarray*}
&&S_{\mathrm{NIE}}^{\mathrm{eff},\mathrm{nonpar}} ( \delta_{1},\theta_{0}
)\\
&&\qquad=s_{\mathrm{NIE}}^{\mathrm{eff},\mathrm{nonpar}} ( O;\delta_{1},\theta_{0} )\\
&&\qquad=\frac{I(E=1)}{f_{E|X}(1|X)} \biggl\{ Y-\eta ( 1,1,X )\\
&&\hspace*{86pt}{} -\frac{%
f_{M|E,X} ( M|E=0,X ) }{f_{M|E,X} ( M|E=1,X ) } \{ Y-%
\mathbb{E} ( Y|X,M,E=1 ) \}\biggr \} \\
&&\quad\qquad{}-\frac{I(E=0)}{f_{E|X}(0|X)} \{ \mathbb{E} ( Y|X,M,E=1 )
-\eta ( 1,0,X ) \}\\
&&\qquad\quad{} +\eta ( 1,1,X ) -\eta (
1,0,X ) +\theta_{0}-\delta_{1}.
\end{eqnarray*}
Thus, the semiparametric efficiency bound for estimating
the natural direct and the natural indirect effects in $\mathcal{M}_{\mathrm
{nonpar}}$ are respectively given by $\mathbb{E} \{ S_{\mathrm{NDE}}^{\mathrm{eff},\mathrm{nonpar}} (
\theta_{0},\delta_{0} ) ^{2} \} ^{-1}$ and $\mathbb{E}%
\{ S_{\mathrm{NIE}}^{\mathrm{eff},\mathrm{nonpar}} ( \delta_{1},\theta_{0} )
^{2} \} ^{-1}$.
\end{theorem}

Although not presented here, Theorem~\ref{teo1} is easily extended to obtain the
efficient influence functions and the respective semiparametric efficiency
bounds for the direct and indirect effects on the risk ratio and the odds
ratio scales by a straightforward application of the delta method. An
important implication of the theorem is that all regular and asymptotically
linear (RAL) estimators of $\theta_{0}$, $\delta_{1}-\theta_{0}$ and
$\theta
_{0}-\delta_{0}$ in model $\mathcal{M}_{\mathrm{nonpar}}$ share the
common influence
functions $S_{\theta_{0}}^{\mathrm{eff},\mathrm{nonpar}} ( \theta_{0} )
,S_{\mathrm{NDE}}^{\mathrm{eff},\mathrm{nonpar}} ( \theta_{0},\delta_{0} ) $ and $%
S_{\mathrm{NIE}}^{\mathrm{eff},\mathrm{nonpar}} ( \delta_{1},\theta_{0} ) $, respectively.
Specifically, any RAL estimator $\widehat{\theta}_{0}$ of the M-functional
$\theta_{0}$ in model $\mathcal{M}_{\mathrm{nonpar}}$, shares a common
asymptotic
expansion,
\[
n^{1/2} ( \widehat{\theta}_{0}-\theta_{0} ) =n^{1/2}\mathbb{P}%
_{n}S_{\theta_{0}}^{\mathrm{eff},\mathrm{nonpar}} ( \theta_{0} ) +o_{P} (
1 ) ,
\]
where $\mathbb{P}_{n} [ \cdot ] =n^{-1}\sum_{i} [ \cdot ]
_{i}$. To illustrate this property of nonparametric RAL estimators, and
as a
motivation to multiply robust estimation when nonparametric methods are not
appropriate, we provide a detailed study of three nonparametric strategies
for estimating the M-functional in a simple yet instructive setting in which
$X$ and $M$ are both discrete with finite support.

\textit{Strategy} 1: The first strategy entails obtaining the maximum
likelihood estimator upon evaluating the M-functional under the empirical
law of the observed data,
\[
\widehat{\theta}_{0}^{\,\mathrm{ym}}= \mathbb{P}_{n}\sum_{m\in\mathcal
{S}}\widehat{%
\mathbb{E}} ( Y|E=1,M=m,X ) \widehat{f}_{M|E,X} (
m|E=0,X ) ,
\]
where $\widehat{f}_{Y|E,M,X}$ and $\widehat{f}_{M|E,X}$ are the empirical
probability mass functions, and $\widehat{\mathbb{E}}(Y|E=e, M=m,X=x)$ is
the expectation of $Y$ under $\widehat{f}_{Y|E,M,X}$.

\textit{Strategy} 2: The second strategy is based on the following
alternative representation of the M-functional:
\begin{eqnarray*}
&&\mathop{\int\hspace*{-4pt}\int}_{\mathcal{S\times X}}
\mathbb{E} ( Y|E=1,M=m,X=x )\, dF_{M|E} ( m|E=0,X=x )
\,dF_{X} ( x )\\
&&\qquad=\sum_{e=0}^{1}\mathop{\int\hspace*{-4pt}\int}_{\mathcal{S\times X}}
\mathbb{E} ( Y|E=1,M=m,X=x ) \frac{I(e=0)}{%
f_{E|X}(e|X=x)}\,dF_{M,E,X} ( m,e,x )\\
&&\qquad=\mathbb{E} \biggl\{ \frac{I(E=0)}{f_{E|X}(0|X)}\mathbb{E} (
Y|E=1,M,X ) \biggr\}.
\end{eqnarray*}

Thus, our second estimator takes the form
\[
\widehat{\theta}_{0}^{\,\mathrm{ye}}=\mathbb{P}_{n} \biggl\{ \frac{I(E=0)}{\widehat{f}%
_{E|X}(0|X)}\widehat{\mathbb{E}} ( Y|E=1,M,X ) \biggr\} ,
\]
with $\widehat{f}_{E|X}$ the empirical estimate of the probability mass
function $f_{E|X}$.

\textit{Strategy} 3: The last strategy is based on a third representation
of the M-functional
\begin{eqnarray*}
&& \mathop{\int\hspace*{-4pt}\int}_{\mathcal{S\times X}}
\mathbb{E} ( Y|E=1,M=m,X=x )\, dF_{M|E} ( m|E=0,X=x )
\,dF_{X} ( x )\\
&&\qquad=\sum_{e=0}^{1}\mathop{\int\hspace*{-4pt}\int\hspace*{-4pt}\int}_{\mathcal{Y\times S\times X}}
y\frac{I(e=1)}{f_{E|X}(e|X=x)}\frac{%
f_{M|E,X} ( M|E=0,X ) }{f_{M|E,X} ( M|E,X ) }%
\,dF_{Y,M,E,X} ( y,m,e,x )\\
&&\qquad=\mathbb{E} \biggl\{ Y\frac{I(E=1)}{f_{E|X}(E|X)}\frac{%
f_{M|E,X} ( M|E=0,X ) }{f_{M|E,X} ( M|E,X ) } \biggr\}.
\end{eqnarray*}

Thus, our third estimator takes the form
\[
\widehat{\theta}_{0}^{\,\mathrm{em}}=\mathbb{P}_{n} \biggl\{ Y\frac{I(E=1)}{\widehat{f}%
_{E|X}(E|X)}\frac{\widehat{f}_{M|E,X} ( M|E=0,X ) }{\widehat{f}%
_{M|E,X} ( M|E,X ) } \biggr\} .
\]
At first glance the three estimators $\widehat{\theta}_{0}^{\,\mathrm{em}}$, $\widehat{\theta}_{0}^{\,\mathrm{ye}}$ and $\widehat{\theta}_{0}^{\,\mathrm{ym}}$ might appear to be
distinct; however, we observe that provided the empirical distribution
function $\widehat{F}_{O}= \widehat{F}_{Y|E,M,X}\times\widehat{F}%
_{M|E,X}\times\widehat{F}_{E|X}\times\widehat{F}_{X}$ satisfies the
positivity assumption, and thus $\widehat{F}_{O}\in\mathcal{M}_{\mathrm
{nonpar}}$,
then actually $\widehat{\theta}_{0}^{\,\mathrm{em}}=\widehat{\theta}_{0}^{\,\mathrm{ye}} =%
\widehat{\theta}_{0}^{\,\mathrm{ym}}=\theta_{0} ( \widehat{F}_{O} ) $ since
the three representations agree on the nonparametric model $\mathcal{M}%
_{\mathrm{nonpar}}$. Therefore we may conclude that these three
estimators are in
fact asymptotically efficient in $\mathcal{M}_{\mathrm{nonpar}}$ with common
influence function $S_{\theta_{0}}^{\mathrm{eff},\mathrm{nonpar}} ( \theta_{0} ) $.
Furthermore, from this observation, one further concludes that (asymptotic)
inferences obtained using one of the three representations are
identical to
inferences using either of the other two representations.

At this juncture, we note that the above equivalence no longer applies when
as we have previously argued will likely occur in practice, $(M,X)$ contains
3 or more continuous variables and/or $X$ is too high dimensional for models
to be saturated or nonparametric, and thus parametric (or
semiparametric)
models are specified for dimension reduction. Specifically, for such
settings, we observe that three distinct modeling strategies are available.
Under the first strategy, the estimator\vspace*{1pt} $\widehat{\theta
}_{0}^{\,\mathrm{ym},\mathrm{par}}$ is
obtained $\widehat{\theta}_{0}^{\,\mathrm{ym}}$ using parametric model
estimates $\widehat{\mathbb{E}}^{\mathrm{par}} ( Y|E,M,X ) $ and $\widehat{f}%
_{M|E,X}^{\,\mathrm{par}} ( m|E,X ) $ instead of their nonparametric
counterparts; similarly under the second strategy,
the estimator $\widehat{\theta}_{0}^{\,\mathrm{ye},\mathrm{par}}$ is obtained\vadjust{\goodbreak}
similarly to
$\widehat{\theta}_{0}^{\,\mathrm{ye}}$
using estimates of parametric models $\widehat{\mathbb{E}}^{\mathrm{par}} (
Y|E=1,M=m,X ) $ and $\widehat{f}_{E|X}^{\,\mathrm{par}}(e|X)$ and finally, under
the third strategy, $\widehat{\theta}_{0}^{\,\mathrm{em},\mathrm{par}}$ is obtained similarly to
$\widehat{\theta}_{0}^{\,\mathrm{em}}$ using $\widehat{f}_{E|X}^{\,\mathrm{par}}(e|X)$ and $%
\widehat{f}_{M|E,X}^{\,\mathrm{par}} ( m|E,X ) $.
Then, it follows that $%
\widehat{\theta}_{0}^{\,\mathrm{ym},\mathrm{par}}$ is CAN under the submodel $\mathcal{M}_{a}$,
but is generally inconsistent if either $\widehat{\mathbb{E}}^{\mathrm{par}} (
Y|E,M,X ) $ or $\widehat{f}_{M|E,X}^{\,\mathrm{par}} ( m|E,X ) $ fails to
be consistent. Similarly, $\widehat{\theta}_{0}^{\,\mathrm{ye},\mathrm{par}}$ and $\widehat
{%
\theta}_{0}^{\,\mathrm{em},\mathrm{par}}$ are, respectively, CAN under the submodels
$\mathcal{M}%
_{b}$ and $\mathcal{M}_{c}$,  but each estimator generally fails to be
consistent outside of the corresponding submodel. In the next section, we
propose an approach that produces a triply robust estimator by
combining the
above three strategies so that only one of models $\mathcal
{M}_{a},\mathcal{M%
}_{b}$ and $\mathcal{M}_{c}$ needs to be valid for consistency of the
estimator.

\subsection{Triply robust estimation}\label{sec2.3}

The proposed triply robust estimator $\widehat{\theta}_{0}^{\,\mathrm{triply}}$
solves%
\[
\mathbb{P}_{n}\widehat{S}_{\theta_{0}}^{\mathrm{eff},\mathrm{nonpar}} ( \widehat{\theta}
_{0}^{\,\mathrm{triply}} ) =0,
\]
where $\widehat{S}_{\theta_{0}}^{\mathrm{eff},\mathrm{nonpar}} ( \theta ) $ is
equal to $S_{\theta_{0}}^{\mathrm{eff},\mathrm{nonpar}} ( \theta ) $ evaluated at $%
\{\widehat{\mathbb{E}}^{\mathrm{par}} ( Y|E,M,X ) $, $\widehat{f}%
_{M|E,X}^{\,\mathrm{par}} ( m|E,X ) $, $\widehat{f}_{E|X}^{\,\mathrm{par}}(e|X)\}$; that
is,%
\begin{eqnarray}\label{triply}
\widehat{\theta}_{0}^{\,\mathrm{triply}} &=&\mathbb{P}_{n} \biggl[ \frac{I \{
E=1 \} \widehat{f}_{M|E,X}^{\,\mathrm{par}} ( M|E=0,X ) }{\widehat{f}%
_{E|X}^{\,\mathrm{par}}(1|X)\widehat{f}_{M|E,X}^{\,\mathrm{par}} ( M|E=1,X ) }\nonumber\\
 &&\hspace*{14pt}{}\times\{ Y-%
\widehat{\mathbb{E}}^{\mathrm{par}} ( Y|X,M,E=1 )  \}
 \nonumber\\[-8pt]\\[-8pt]
&&\hspace*{14pt}{} +\frac{I(E=0)}{\widehat{f}_{E|X}^{\,\mathrm{par}}(0|X)} \{ \widehat{%
\mathbb{E}}^{\mathrm{par}} ( Y|X,M,E=1 )\nonumber\\
&&\hspace*{32pt}{} -\widehat{\eta}^{\,\mathrm{par}} (
1,0,X ) \} +\widehat{\eta}^{\,\mathrm{par}} ( 1,0,X ) \biggr] ,
\nonumber
\end{eqnarray}
  is CAN in model $\mathcal{M}_{\mathrm{union}}=\mathcal
{M}_{a}\cup
\mathcal{M}_{b}\cup\mathcal{M}_{c}$, where
\[
\widehat{\eta}^{\,\mathrm{par}} ( e,e^{\ast},X ) =\int_{\mathcal{S}}\widehat{%
\mathbb{E}}^{\mathrm{par}} ( Y|X,M=m,E=e ) \widehat{f}_{M|E,X}^{\,\mathrm{par}} (
m|E=e^{\ast},X )\, d\mu(m).
\]

In the next theorem, the estimator in the above display is combined
with a
doubly robust estimator $\widehat{\delta}_{e}^{\,\mathrm{doubly}}$ of $\delta_{e}$
[see \citet{vanderLaan} or \citet{Tsiatis}], to obtain
multiply robust estimators of natural direct and indirect effects, where
\[
\widehat{\delta}_{e}^{\,\mathrm{doubly}}=\mathbb{P}_{n} \biggl[ \frac{I(E=e)}{\widehat
{f}%
_{E|X}^{\,\mathrm{par}}(e|X)} \{ Y-\widehat{\eta}^{\,\mathrm{par}} ( e,e,X )
\} +\widehat{\eta}^{\,\mathrm{par}} ( e,e,X ) \biggr] .
\]
To state the result, we set $\widehat{\mathbb{E}}^{\mathrm{par}} ( Y|X,M,E )
=\mathbb{E}^{\mathrm{par}} ( Y|X,M,E;\widehat{\beta}_{y} ) =\break g^{-1} (
\widehat{\beta}_{y}^{T}h(X, M,E) ) $, where\vadjust{\goodbreak} $g$ is a known link
function, and $h$ is a user specified function of ($X,M,E)$ so that
$\mathbb{E}%
^{\mathrm{par}} ( Y|X,M,E;\beta_{y} ) =g^{-1} ( \beta
_{y}^{T}h(X,M,E) ) $ entails a working regression model for $\mathbb{E}%
( Y|X,M,E ) $, and $\widehat{\beta}_{y}$ solves the estimating
equation
\[
0=\mathbb{P}_{n} [ S_{y} ( \widehat{\beta}_{y} ) ] =%
\mathbb{P}_{n} \bigl[ h(X,M,E) \bigl( Y-g^{-1} ( \widehat{\beta}%
_{y}^{T}h(X,M,E) ) \bigr) \bigr] .
\]
Similarly, we set $\widehat{f}_{M|E,X}^{\,\mathrm{par}} ( m|E,X )
=f_{M|E,X}^{\,\mathrm{par}} ( m|E,X;\widehat{\beta}_{m} ) $ for $%
f_{M|E,X}^{\,\mathrm{par}} ( m|E,X;\break \beta_{m} ) $, a parametric model for the
density of $[M|E,X]$ with $\widehat{\beta}_{m}$, solving
\[
0=\mathbb{P}_{n} [ S_{m} ( \widehat{\beta}_{m} ) ] =%
\mathbb{P}_{n} \biggl[ \frac{\partial}{\partial\beta_{m}}\log
f_{M|E,X}^{\,\mathrm{par}} ( M|E,X;\widehat{\beta}_{m} ) \biggr] ,
\]
and we set $\widehat{f}_{E|X}^{\,\mathrm{par}}(e|X)=f_{E|X}^{\,\mathrm{par}}(e|X;\widehat
{\beta}%
_{e})$ for $f_{E|X}^{\,\mathrm{par}}(e|X;\beta_{e})$, a parametric model for the
density of $[E|X]$ with $\widehat{\beta}_{e}$ solving
\[
0=\mathbb{P}_{n} [ S_{e} ( \widehat{\beta}_{e} ) ] =%
\mathbb{P}_{n} \biggl[ \frac{\partial}{\partial\beta_{e}}\log
f_{E|X}^{\,\mathrm{par}} ( E|X;\widehat{\beta}_{e} ) \biggr] .
\]

\begin{theorem}\label{teo2}
Suppose that the assumptions of Theorem~\ref{teo1} hold, and that
the regularity conditions stated in the \hyperref[app]{Appendix}
hold and that $\beta
_{m},\beta_{e}$ and $\beta_{y}$ are variation independent.
\begin{longlist}[(iii)]
\item[(i)] Mediation functional: Then, $\sqrt{n}(\widehat{\theta}%
_{0}^{\,\mathrm{triply}}-\theta_{0})$ is RAL under model $\mathcal{M}%
_{\mathrm{union}} $ with influence function %
\begin{eqnarray*}
&&S_{\theta_{0}}^{\mathrm{union}} ( \theta_{0},\beta^{\ast} ) \\
&&\qquad=S_{\theta_{0}}^{\mathrm{eff},\mathrm{nonpar}} ( \theta_{0},\beta^{\ast} ) -%
\frac{\partial\mathbb{E} \{ S_{\theta_{0}}^{\mathrm{eff},\mathrm{nonpar}} ( \theta
_{0},\beta ) \} }{\partial\beta^{T}}\bigg|_{\beta^{\ast}}\mathbb{E%
} \biggl\{ \frac{\partial S_{\beta} ( \beta ) }{\partial\beta^{T}}%
\bigg|_{\beta^{\ast}} \biggr\} ^{-1}S_{\beta} ( \beta^{\ast} ) ,
\end{eqnarray*}
and thus converges in distribution to a $N ( 0,\Sigma_{\theta
_{0}} ) $, where
\[
\Sigma_{\theta_{0}} ( \theta_{0},\beta^{\ast} ) =\mathbb{E}%
( S_{\theta_{0}}^{\mathrm{union}} ( \theta_{0},\beta^{\ast} )
^{2} ) ,
\]
with $\beta^{T}= ( \beta_{m}^{T},\beta_{e}^{T},\beta
_{y}^{T} ) $ and $S_{\beta} ( \beta ) = (
S_{m}^{T} ( \beta_{m} ) ,S_{e}^{T} ( \beta_{e} )
,S_{y}^{T} ( \beta_{y} ) ) ^{T}$, and with~$\beta
^{\ast}$ denoting the probability limit of the estimator $%
\widehat{\beta}= ( \widehat{\beta}_{m}^{T},\widehat{\beta}_{e}^{T},%
\widehat{\beta}_{y}^{T} ) ^{T}$.

\item[(ii)]Natural direct effect: Similarly, $\sqrt{n}(\widehat{\theta
}%
_{0}^{\,\mathrm{triply}}-\widehat{\delta}_{0}^{\,\mathrm{doubly}}- ( \theta_{0}-\delta
_{0} ) )$ is RAL under model $\mathcal{M}_{\mathrm{union}}$
with influence function $S_{\mathrm{NDE}}^{\mathrm{union}} ( \theta_{0},\delta
_{0},\beta^{\ast} ) $ defined as $S_{\theta
_{0}}^{\mathrm{union}} ( \theta_{0},\beta^{\ast} ) $ with $%
S_{\mathrm{NDE}}^{\mathrm{eff},\mathrm{nonpar}} ( \theta_{0},\delta_{0},\beta^{\ast} ) $
 replacing $S_{\theta_{0}}^{\mathrm{eff},\mathrm{nonpar}} ( \theta_{0},\beta
^{\ast} )$,   and asymptotic variance $\Sigma_{\theta
_{0}-\delta_{0}} ( \delta_{1},\theta_{0},\beta^{\ast} ) $
defined accordingly.

\item[(iii)]Natural indirect effect: Similarly, $\sqrt{n}(\widehat
{\delta}%
_{1}^{\,\mathrm{doubly}}-\widehat{\theta}_{0}^{\,\mathrm{triply}}- ( \delta_{1}-\theta
_{0} ) )$ is RAL under model $M_{\mathrm{union}}$ with
influence function $S_{\mathrm{NIE}}^{\mathrm{union}} ( \delta_{1},\theta_{0},\beta
^{\ast} ) $ defined as $S_{\theta_{0}}^{\mathrm{union}} ( \theta
_{0},\beta^{\ast} ) $ with $S_{\mathrm{NIE}}^{\mathrm{eff},\mathrm{nonpar}} (
\delta_{1},\theta_{0},\beta^{\ast} ) $ replacing $%
S_{\theta_{0}}^{\mathrm{eff},\mathrm{nonpar}} ( \theta_{0},\beta^{\ast} ) $, and asymptotic variance $\Sigma_{\delta_{1}-\theta
_{0}} ( \delta_{1},\theta_{0},\beta^{\ast} ) $ defined
accordingly.

\item[(iv)]$\widehat{\theta}_{0}^{\,\mathrm{triply}}$, $\widehat
{\theta}%
_{0}^{\,\mathrm{triply}}-\widehat{\delta}_{0}^{\,\mathrm{doubly}}$ and $\widehat
{\delta
}_{1}^{\,\mathrm{doubly}}-\widehat{\theta}_{0}^{\,\mathrm{triply}}$ are semiparametric
locally efficient in the sense that they are RAL under model $\mathcal
{M}%
_{\mathrm{union}}$ and respectively achieve the semiparametric efficiency\vadjust{\goodbreak}
bound for $\theta_{0}$, $\theta_{0}-\delta_{0}$, and $\delta
_{1}-\theta
_{0}$ under model $\mathcal{M}_{\mathrm{union}}$ at the
intersection submodel $\mathcal{M}_{a}\cap\mathcal{M}_{b}\cap
\mathcal{M}_{c}$, with respective efficient influence
functions: $S_{\theta_{0}}^{\mathrm{eff},\mathrm{nonpar}} ( \theta_{0},\beta^{\ast
} ) ,S_{\mathrm{NDE}}^{\mathrm{eff},\mathrm{nonpar}} ( \theta_{0},\delta_{0},\beta^{\ast
} ) $ and $S_{\mathrm{NIE}}^{\mathrm{eff},\mathrm{nonpar}} ( \delta_{1},\theta
_{0},\beta^{\ast} ) $.
\end{longlist}
\end{theorem}

Empirical versions of $\Sigma_{\theta_{0}-\delta_{0}} ( \delta
_{1},\theta_{0},\beta^{\ast} ) $ and $\Sigma_{\delta_{1}-\theta
_{0}} ( \delta_{1},\theta_{0},\beta^{\ast} ) $ are easily
obtained, and the corresponding Wald-type confidence intervals can be used
to make formal inferences about natural direct and indirect effects. It is
also straightforward to extend the approach to the risk ratio and odds ratio
scales for binary $Y$. By a theorem due to \citet{Robins8}, part
(iv) of the theorem implies that when all models are correct, $\widehat{
\theta}_{0}^{\,\mathrm{triply}}$, $\widehat{\theta
}_{0}^{\,\mathrm{triply}}-\widehat{%
\delta}_{0}^{\,\mathrm{doubly}}$ and $\widehat{\delta}_{1}^{\,\mathrm{doubly}}-%
\widehat{\theta}_{0}^{\,\mathrm{triply}}$ are semiparametric efficient in model $%
\mathcal{M}_{\mathrm{nonpar}}$ at the intersection submodel $\mathcal{M}
_{a}\cap\mathcal{M}_{b} \cap \mathcal{M}_{c}$.

\section{A simulation study of estimators of direct effect}\label{sec3}

In this section, we report a simulation study which illustrates the finite
sample performance of the various estimators described in previous sections.
We generated 1000 samples of size $n=600,1000$ from the following model:
\begin{enumerate}[(Model.M)]
\item[(Model.X)] $X_{1}\sim \mathit{Bernoulli}(0.4);[X_{2}|X_{1}]\sim
\mathit{Bernoulli}(0.3+0.4X_{1})$;

\item[] $ \lbrack X_{3}|X_{1},X_{2}]\sim
-0.024-0.4X_{1}+0.4X_{2}+N(0,1)$;

\item[(Model.E)] $[E|X_{1},X_{2},X_{3}]\sim \mathit{Bernoulli} ( [1+\exp\{- (
0.4+X_{1}-X_{2}+0.1X_{3}-1.5X_{1}X_{3} ) \}]^{-1} ) $;

\item[(Model.M)] $[M|E,X_{1},X_{2},X_{3}]\sim
\mathit{Bernoulli}([1+\exp\{-(0.5-X_{1}+0.5X_{2}$

\item[] $ -0.9X_{3}+E-1.5X_{1}X_{3} )\}]^{-1});$

\item[(Model.Y)] $[Y|M,E,X_{1},X_{2},X_{3}]\sim1+0.2X_{1}+0.3X_{2}+1.4X_{3}$

\item[] $ -2.5E-3.5M+5EM+N(0,1)$.
\end{enumerate}

We then evaluated the performance of the following four estimators of the
natural direct effect $\widehat{\theta}_{0}^{\,\mathrm{em}}-\widehat{\delta}%
_{0}^{\,\mathrm{doubly}}$, $\widehat{\theta}_{0}^{\,\mathrm{ye}}-\widehat{\delta
}_{0}^{\,\mathrm{doubly}}$, $%
\widehat{\theta}_{0}^{\,\mathrm{ym}}-\widehat{\delta}_{0}^{\,\mathrm{doubly}}$ and $\widehat
{%
\theta}_{0}^{\,\mathrm{triply}}-\widehat{\delta}_{0}^{\,\mathrm{doubly}}$. Note that the doubly
robust estimator $\widehat{\delta}_{0}^{\,\mathrm{doubly}}$ was used throughout to
estimate $\delta_{0}=\mathbb{E} ( Y_{0} ) $. To assess the impact
of modeling error, we evaluated these estimators in four separate scenarios.
In the first scenario, all models were correctly specified, whereas the
remaining three scenarios respectively mis-specified only one of Model E,
Model M and Model Y. In order to mis-specify Model E and Model~M, we
respectively left out the $X_{1}X_{3}$ interaction when fitting each model,
and we assumed an incorrect log--log link function. The incorrect model
for $Y $ simply assumed no $EM$ interaction.

\begin{table}
\caption{Simulation results $n=600$}\label{tab1}
\begin{tabular*}{\textwidth}{@{\extracolsep{\fill}}lld{2.3}d{2.3}d{2.4}c@{}}
\hline
& & \multicolumn{1}{c}{$\bolds{\mathcal{M}_{\mathrm{ym}}}$} & \multicolumn{1}{c}{$\bolds{\mathcal{M}_{\mathrm{ye}}}$} & \multicolumn{1}{c}{$\bolds{\mathcal{M}_{\mathrm{em}}}$} & $\bolds{\mathcal{M}_{\mathrm{union}}}$ \\
\hline
All correct & bias & 0.002 & 0.008 & 0.002 &0.005 \\
$ $ & $MC$ s.e.\tabnoteref[\ast]{vz} & 0.005 & 0.007 & 0.006 & 0.006 \\
$Y$ wrong & bias & -0.500 & -0.500 & 0.0001 & 0.004 \\
$ $ & $MC$ s.e. & 0.005 & 0.006 & 0.006 & 0.006 \\
$M$ wrong & bias & 0.038 & 0.008 & -0.054 &0.003 \\
& $MC$ s.e. & 0.005 & 0.007 & 0.006 & 0.006 \\
$E$ wrong & bias & 0.003 & 0.027 & 0.059 &0.004 \\
$ $ & $MC$ s.e. & 0.005 & 0.005 & 0.005 & 0.005 \\
\hline
\end{tabular*}
\tabnotetext[]{tz}{$\mathcal{M}_{\mathrm{ym}}$: $\widehat{\theta}_{0}^{\,\mathrm{ym}}-\widehat{\delta}%
_{0}^{\,\mathrm{doubly}}$; $\mathcal{M}_{\mathrm{ye}}$: $\widehat{\theta}_{0}^{\,\mathrm{ye}}-\widehat
{\delta}_{0}^{\,\mathrm{doubly}}$; $\mathcal{M}_{\mathrm{em}}$: $\widehat{\theta}_{0}^{\,\mathrm{em}}-\widehat{%
\delta}_{0}^{\,\mathrm{doubly}}$;  $\mathcal{M}_{\mathrm{union}}$: $\widehat{\theta
}_{0}^{\,\mathrm{triply}}-\widehat{\delta}_{0}^{\,\mathrm{doubly}}$.}
\tabnotetext[\ast]{vz}{Monte Carlo standard error.}\vspace*{-2pt}
\end{table}

\begin{table}[b]\vspace*{-2pt}
\caption{Simulation results $n=1000$}\label{tab2}
\begin{tabular*}{\textwidth}{@{\extracolsep{\fill}}lld{2.3}d{2.3}d{2.3}c@{}}
\hline
& & \multicolumn{1}{c}{$\bolds{\mathcal{M}_{\mathrm{ym}}}$} & \multicolumn{1}{c}{$\bolds{\mathcal{M}_{\mathrm{ye}}}$} & \multicolumn{1}{c}{$\bolds{\mathcal{M}_{\mathrm{em}}}$} & $\bolds{\mathcal{M}_{\mathrm{union}}}$ \\
\hline
All correct & bias & 0.001 & 0.009 &   0.001 &0.001 \\
$ $ & $MC$ s.e.\tabnoteref[\ast]{vz} & 0.004 & 0.005 & 0.004 & 0.004 \\
$Y$ wrong & bias & -0.484 & -0.484 &0.003 & 0.003 \\
$ $ & $MC$ s.e. & 0.004 & 0.004 & 0.004 & 0.004 \\
$M$ wrong & bias & 0.136 & -0.008 & 0.056 &  0.01\phantom{0} \\
& $MC$ s.e. & 0.004 & 0.05 & 0.004 & 0.01\phantom{0} \\
$E$ wrong & bias & 0.001 & -0.024 &   -0.054& 0.001 \\
$ $ & $MC$ s.e. & 0.004 & 0.004 & 0.004 & 0.004 \\
\hline
\end{tabular*}
\tabnotetext[]{tz}{$\mathcal{M}_{\mathrm{ym}}$: $\widehat{\theta}_{0}^{\,\mathrm{ym}}-\widehat{\delta}%
_{0}^{\,\mathrm{doubly}}$; $\mathcal{M}_{\mathrm{ye}}$: $\widehat{\theta}_{0}^{\,\mathrm{ye}}-\widehat
{\delta
}_{0}^{\,\mathrm{doubly}}$; $\mathcal{M}_{\mathrm{em}}$: $\widehat{\theta}_{0}^{\,\mathrm{em}}-\widehat{%
\delta}_{0}^{\,\mathrm{doubly}}$; $\mathcal{M}_{\mathrm{union}}$ : $\widehat{\theta
}_{0}^{\,\mathrm{triply}}-\widehat{\delta}_{0}^{\,\mathrm{doubly}}$.}
\tabnotetext[\ast]{vz}{Monte Carlo standard error.}
\end{table}

Tables~\ref{tab1} and~\ref{tab2} summarize the simulation results which largely agree
with the
theory developed in the previous sections. Mainly, all proposed estimators
performed well at both moderate and large sample sizes in the absence of
modeling error. Furthermore, under the partially mis-specified model in
which Model.Y was incorrect, both estimators, $\widehat{\theta
}_{0}^{\,\mathrm{ye}}-%
\widehat{\delta}_{0}^{\,\mathrm{doubly}}$ and $\widehat{\theta}_{0}^{\,\mathrm{ym}}-\widehat
{%
\delta}_{0}^{\,\mathrm{doubly}}$, showed significant bias irrespective of sample size,
while $\widehat{\theta}_{0}^{\,\mathrm{em}}-\widehat{\delta}_{0}^{\,\mathrm{doubly}}$ and\vadjust{\goodbreak} $%
\widehat{\theta}_{0}^{\,\mathrm{triply}}-\widehat{\delta}_{0}^{\,\mathrm{doubly}}$ both
performed well. Similarly when Model~M was incorrect, the estimators $%
\widehat{\theta}_{0}^{\,\mathrm{em}}-\widehat{\delta}_{0}^{\,\mathrm{doubly}}$ and $\widehat
{%
\theta}_{0}^{\,\mathrm{ym}}-\widehat{\delta}_{0}^{\,\mathrm{doubly}}$ resulted in large bias,
when compared to the relatively small bias of $\widehat{\theta
}_{0}^{\,\mathrm{ye}}-%
\widehat{\delta}_{0}^{\,\mathrm{doubly}}$ and $\widehat{\theta
}_{0}^{\,\mathrm{triply}}-\widehat{%
\delta}_{0}^{\,\mathrm{doubly}}$. Finally, mis-specifying Model E lead to
estimators $%
\widehat{\theta}_{0}^{\,\mathrm{ye}}-\widehat{\delta}_{0}^{\,\mathrm{doubly}}$ and $\widehat
{%
\theta}_{0}^{\,\mathrm{em}}-\widehat{\delta}_{0}^{\,\mathrm{doubly}}$ that were significantly
more biased than the estimators $\widehat{\theta}_{0}^{\,\mathrm{ym}}-\widehat
{\delta}%
_{0}^{\,\mathrm{doubly}}$ and $\widehat{\theta}_{0}^{\,\mathrm{triply}}-\widehat{\delta}%
_{0}^{\,\mathrm{doubly}}$. Interestingly, the efficiency loss of the multiply robust
estimator remained relatively small when compared to the consistent
nonrobust estimator under the various scenarios, suggesting that, at least
in this simulation study, the benefits of robustness appear to outweigh the
loss of efficiency.

\section{A data application}\label{sec4}

In this section, we illustrate the methods in a~real world application from
the psychology literature on mediation. We re-analyze data from The Job
Search Intervention Study (JOBS II) also analyzed by \citet{Imai2}.\vadjust{\goodbreak}
JOBS II is a randomized field experiment that investigates the efficacy
of a
job training intervention on unemployed workers. The program is
designed not
only to increase reemployment among the unemployed but also to enhance the
mental health of the job seekers. In the study, 1801 unemployed workers
received a pre-screening questionnaire and were then randomly assigned to
treatment and control groups. The treatment group with $E=1$
participated in
job skills workshops in which participants learned job search skills and
coping strategies for dealing with setbacks in the job search process. The
control group with $E=0$ received a booklet describing job search tips. An
analysis considers a continuous outcome measure $Y $of depressive
symptoms based on the Hopkins Symptom Checklist [\citet{Imai2}]. In the JOBS II
data, a
continuous measure of job search self-efficacy represented the hypothesized
mediating variable~$M$. The data also included baseline covariates $X$
measured before administering the treatment including: pretreatment
level of
depression, education, income, race, marital status, age, sex, previous
occupation, and the level of economic hardship.

Note that by randomization, the density of $[E|X]$ was known by design not
to depend on covariates, and therefore its estimation is not prone to
modeling error. The continuous outcome and mediator variables were modeled
using linear regression models with Gaussian error, with main effects
for $(E,M,X)$ included in the outcome regression and main effects for $(E,X)$
included in the mediator regression. Table~\ref{tab3} summarizes results obtained
using $\widehat{\theta}_{0}^{\,\mathrm{em}}$, $\widehat{\theta}_{0}^{\,\mathrm{ye}}$, $\widehat
{\theta}_{0}^{\,\mathrm{ym}}$ and $\widehat{\theta}_{0}^{\,\mathrm{triply}}$ together with $%
\widehat{\delta}_{e}^{\,\mathrm{doubly}}$, $e=0,1$, to estimate the direct and
indirect effects of the treatment.

\begin{table}
\caption{Estimated causal effects of interest using the job search
intervention study data}\label{tab3}
\begin{tabular*}{\textwidth}{@{\extracolsep{\fill}}lld{2.4}d{2.4}d{2.4}d{2.4}@{}}
\hline
& & \multicolumn{1}{c}{$\bolds{\mathcal{M}_{\mathrm{ym}}}$} & \multicolumn{1}{c}{$\bolds{\mathcal{M}_{\mathrm{ye}}}$} & \multicolumn{1}{c}{$\bolds{\mathcal{M}_{\mathrm{em}}}$} & \multicolumn{1}{c@{}}{$\bolds{\mathcal{M}_{\mathrm{union}}}$} \\
\hline
Direct effect & Estimate & -0.0310 & -0.0310 &   0.0280& -0.0409 \\
$ $ & s.e.\tabnoteref[\ast]{vz} &   0.0124 &   0.0620 &  0.0465 & 0.021 7 \\
Indirect effect & Estimate & -0.0160 & -0.0160 & -0.0750 &  -0.0070 \\
& s.e.\tabnoteref[\ast]{vz} &   0.0372 &   0.0620 &0.0434 &   0.021 7 \\
\hline
\end{tabular*}
\tabnotetext[\ast]{vz}{Nonparametric bootstrap standard errors.}
\end{table}

Point estimates of both natural direct and indirect effects closely agreed
under models $\mathcal{M}_{\mathrm{ym}}$ and $\mathcal{M}_{\mathrm{ye}}$, and also agreed with
the results of \citet{Imai2}. We should note that inferences
under our
choice of $\mathcal{M}_{\mathrm{ym}}$ are actually robust to the normality assumption
and, as in \citet{Imai2}, only require that the mean structure of
$[Y|E,M,X]$ and $[ M|E,X]$ is correct. In contrast,
inferences under model $\mathcal{M}_{\mathrm{em}}$ require a~correct model for the
mediator density. This distinction may partly explain the apparent
disagreement in the estimated direct effect under $\mathcal{M}_{\mathrm{em}}$ when
compared to the other methods,\vadjust{\goodbreak} also suggesting that the Gaussian error model
for $M$ is not entirely appropriate. The multiply robust estimate of the
natural direct effect is consistent with estimates obtained under
models $\mathcal{M}_{\mathrm{ym}}$ and $\mathcal{M}_{\mathrm{ye}}$, and is statistically significant,
suggesting that the intervention may have beneficial direct effects on
participants' mental health; while the multiply robust approach
suggests a
much smaller indirect effect than all other estimators although none
achieved statistical significance.

\section{\texorpdfstring{Improving the stability of $\widehat{ \theta}_{0}^{\,\mathrm{triply}}$ when weights are highly variable}
{Improving the stability of theta 0 triply when weights are highly variable}}\label{sec5}

The triply robust estimator $\widehat{\theta}_{0}^{\,\mathrm{triply}}$ which involves
inverse probability weights for the exposure and mediator variables, clearly
relies on the positivity assumption, for good finite sample
performance. But
as recently shown by \citet{Kang} in the context of missing
outcome data, a practical violation of positivity in data analysis can
severely compromise inferences based on such methodology; although their
analysis did not directly concern the M-functional $\theta_{0} $. Thus,
it is crucial to critically examine, as we do below in a simulation study,
the extent to which the various estimators discussed in this paper are
susceptible to a practical violation of the positivity assumption, and to
consider possible approaches to improve the finite sample performance of
these estimators in the context of highly variable empirical weights.
Methodology to enhance the finite sample behavior of $\widehat{\delta}%
_{j}^{\,\mathrm{doubly}}$ is well studied in the literature and is not considered here;
see, for example, \citet{Robins9}, \citet{Cao} and Tan
(\citeyear{T10}). We
first describe an approach to enhance the finite sample performance of $
\widehat{\theta}_{0}^{\,\mathrm{triply}}$, particularly in the presence of highly
variable empirical weights. To focus the exposition, we only consider the
case of a continuous $Y$ and a binary~$M$, but in principle, the approach
could be generalized to a more general setting. The proposed enhancement
involves two modifications.

The first modification adapts to the mediation context, an approach
developed for the missing data context (and for the estimation of total
effects) in \citet{Robins9}. The basic guiding principle of the approach
is to carefully modify the estimation of the outcome and mediator
models in
order to ensure that the triply robust estimator given by equation (\ref{triply}) has the simple M-functional representation
\[
\widehat{\theta}_{0}^{\,\mathrm{triply},\dag}=\mathbb{P}_{n} \{ \widehat{\eta}%
^{\mathrm{par},\dag} ( 1,0,X ) \},
\]
where $\widehat{\eta}^{\,\mathrm{par},\dag} ( 1,0,X ) $ is carefully
estimated to ensure multiple robustness. The reason for favoring an
estimator with the above representation is that it is expected to be more
robust to practical positivity violation because it does not directly depend
on inverse probability weights. However, as we show next, to ensure multiple
robustness, estimation of $\eta^{\mathrm{par}}$ involves inverse probability
weights, and therefore, $\widehat{\theta}_{0}^{\,\mathrm{triply},\dag}$ indirectly
depends on such weights. Our strategy involves a second step to
minimize the
potential impact of this indirect dependence on weights.\vadjust{\goodbreak}

In the following, we assume, to simplify the exposition, that a simple linear
model is used:
\[
\mathbb{E}^{\mathrm{par}} ( Y|X,M,E=1 ) =\mathbb{E}^{\mathrm{par}} (
Y|X,M,1;\beta_{y} ) = [ 1,X^{T},M ] \beta_{y}.
\]
Then, similar to \citet{Robins9}, one can verify that the above
M-functional representation of a triply robust estimator is obtained by
estimating $f_{M|E,X}^{\,\mathrm{par}} ( M|E=0,X ) $ with $\widehat{f}%
_{M|E,X}^{\,\mathrm{par},\dag} ( M|E=0,X ) $ obtained via weighted logistic
regression in the unexposed-only, with weight $\widehat{f}%
_{E|X}^{\,\mathrm{par}}(0|X)^{-1}$; and by estimating $\mathbb{E}^{\mathrm{par}} (
Y|X,M,E=1 ) $ using weighted OLS of $Y$ on $(M,X)$ in the exposed-only,
with weight
\[
\widehat{f}_{M|E,X}^{\,\mathrm{par},\dag} ( M|E=0,X ) \{\widehat{f}%
_{E|X}^{\,\mathrm{par}}(1|X)\widehat{f}_{M|E,X}^{\,\mathrm{par},\dag} ( M|E=1,X )
\}^{-1};
\]
provided that both working models include an intercept. The second
enhancement to minimize undue influence of variable weights on the
M-func\-tional estimator, entails using $\widehat{f}_{E|X}^{\,\mathrm{par},\dag}$
in the
previous step instead of $\widehat{f}_{E|X}^{\,\mathrm{par}}$, where
\[
\operatorname{logit}\widehat{f}_{E|X}^{\,\mathrm{par},\dag}(1|X)=\operatorname{logit}\widehat{f}%
_{E|X}^{\,\mathrm{par}}(1|X)+\widehat{C}_{1}
\]
with
\[
\widehat{C}_{1}=-\log\bigl(1-\mathbb{P}_{n}(E)\bigr)+\log\bigl(\mathbb
{P}_{n}[E\widehat{f}_{E|X}^{\,\mathrm{par}}(0|X)/\widehat{f}_{E|X}^{\,\mathrm{par}}(1|X)]\bigr).
\]
This second modification ensures a certain boundedness property of inverse
propensity score-weighting. Specifically, for any bounded function $R=r(Y,M)$
of $Y$ and $M$; consider for a moment the goal of estimating the
counterfactual mean $\mathbb{E} \{ r(Y_{1},M_{1}) \} $; then it is
well known that even though $R$ is bounded, the simple inverse-probability
weighting estimator $\mathbb{P}_{n} \{ ER\widehat{f}%
_{E|X}^{\,\mathrm{par}}(1|\break X)^{-1} \} $ could easily be unbounded, particularly if
positivity is practically violated. In contrast, as we show next, the
estimator $\mathbb{P}_{n} \{ ER\widehat{f}_{E|X}^{\,\mathrm{par},\dag
}(1|X)^{-1} \} $ is generally bounded. To see why, note that
\begin{eqnarray*}
\mathbb{P}_{n} \{ ER\widehat{f}_{E|X}^{\,\mathrm{par},\dag}(1|X)^{-1} \}
&=&
\mathbb{P}_{n} \{ ER\widehat{f}_{E|X}^{\,\mathrm{par},\dag}(0|X)\widehat{f}%
_{E|X}^{\,\mathrm{par},\dag}(1|X)^{-1} \} +\mathbb{P}_{n} \{ R \} \\
&=&\mathbb{P}_{n} \biggl\{ R\frac{E\widehat{f}_{E|X}^{\,\mathrm{par}}(0|X)\widehat{f}%
_{E|X}^{\,\mathrm{par}}(1|X)^{-1}}{\mathbb{P}_{n}[E\widehat{f}_{E|X}^{\,\mathrm{par}}(0|X)%
\widehat{f}_{E|X}^{\,\mathrm{par}}(1|X)^{-1}]}\bigl(1-\mathbb{P}_{n}(E)\bigr) \biggr\}\\
&&{} +\mathbb{P}%
_{n} \{ R \}
\end{eqnarray*}
which is bounded since the second term is bounded, and the first term
is a
convex combination of bounded variables, and therefore is also bounded.
Furthermore, $\mathbb{P}_{n}[E\widehat{f}_{E|X}^{\,\mathrm{par},\dag}(0|X)\widehat
{f}%
_{E|X}^{\,\mathrm{par},\dag}(1|X)^{-1}]$ converges in probability to $(1-\mathbb{E}(E))
$ provided that $\widehat{f}_{E|X}^{\,\mathrm{par}}$ converges to $f_{E|X}$, ensuring
that the expression in the above display is consistent for $\mathbb{E}%
\{ r(Y_{1},M_{1}) \} $. The nonparametric bootstrap is most
convenient for inference using $\widehat{f}_{E|X}^{\,\mathrm{par},\dag}$.\vadjust{\goodbreak}

In the next section, we study, in the context of highly variable
weights, the
behavior of our previous estimators of $\theta_{0}$, together with
that of
the enhanced estimators $\widehat{\theta}_{0}^{\,\mathrm{triply},\dag,j}=\mathbb
{P}%
_{n} \{ \widehat{\eta}^{\mathrm{par},\dag,j} ( 1,0,X ) \} $, $%
j=1,2$, where $\widehat{\eta}^{\,\mathrm{par},\dag,1}$ is constructed as described
above using $\widehat{f}_{E|X}^{\,\mathrm{par}} $, and $\widehat{\eta}^{\,\mathrm{par},\dag,2}
$ uses $\widehat{f}_{E|X}^{\,\mathrm{par},\dag}$.

\section{A simulation study where positivity is practically violated}\label{sec6}

We adapt\-ed to the mediation setting, the missing data simulation scenarios
in \citet{Kang} which were specifically designed so that, when
misspecified, working models are nonetheless nearly correct, but yield highly
variable inverse probability weights with practical positivity
violation in
the context of estimation. We generated 1000 samples of size $n=200,1000$
from the following model:
\begin{enumerate}[(Model.M)]
\item[(Model.X)] $Z=Z_{1},Z_{2},Z_{3},Z_{4}\stackrel{\mathrm{i.i.d.}}{\sim
}N(0,1);X_{1}=\exp(Z_{1}/2);$

\item[]$ X_{2}=Z_{2}/ \{ 1+\exp(Z_{1}) \} +10$;
$X_{3}=(Z_{1}Z_{3}/25+0.6)^{3}$

\item[] and $X_{4}=(Z_{2}+Z_{4}+20)^{2}$, so that $Z$ may be
expressed in terms of $X$.

\item[(Model.E)] $[E|X_{1},X_{2},X_{3}]\sim \mathit{Bernoulli} ( [1+\exp\{ (
Z_{1}-0.5Z_{2}+0.25Z_{3}$
\item[]$+0.1Z_{4} ) \}]^{-1} ) $;

\item[(Model.M)] $[M|E,X_{1},X_{2},X_{3}]\sim
\mathit{Bernoulli}([1+\exp\{-(0.5-Z_{1}+0.5Z_{2}$

\item[]$ -0.9Z_{3}+Z_{4}-1.5E )\}]^{-1})$;

\item[(Model.Y)] $[Y|M,E,X_{1},X_{2},X_{3}]\sim210+27.4Z_{1}+13.7Z_{3}+13.7Z_{3}$

\item[]$ +M+E+N(0,1)$.
\end{enumerate}

Correctly specified working models were thus achieved when an additive
linear regression of $Y$ on $Z$, a logistic regression of $M$ with linear
predictor additive in $Z$ and $E$ and a logistic regression of $E$ with
linear predictor additive in the $Z$, respectively. Incorrect specification
involved fitting these models with $X$ replacing~$Z$, which produces higly
variable weights. For instance, an estimated propensity score as small as
$5.5 \times10^{-33}$ occurred in the simulation study reflecting an
effective violation of positivity; similarly, a mediator predicted
probability as small as $3\times10^{-20}$ also occured in the simulation
study.

\begin{table}[t!]
\tabcolsep=0pt
\caption{Simulation results $n=200$}\label{tab4}
\begin{tabular*}{\textwidth}{@{\extracolsep{\fill}}lld{2.3}d{3.3}d{4.3}d{3.4}d{2.4}d{2.8}@{}}
\hline
& & \multicolumn{1}{c}{$\bolds{\mathcal{M}_{\mathrm{ym}}}$} & \multicolumn{1}{c}{$\bolds{\mathcal{M}_{\mathrm{ye}}}$} & \multicolumn{1}{c}{$\bolds{\mathcal{M}_{\mathrm{em}}}$} &
\multicolumn{1}{c}{$\bolds{\mathcal{M}_{\mathrm{union}}}$} & \multicolumn{1}{c}{$\bolds{\mathcal{M}_{\mathrm{union}}^{\dag,1}}$} &
\multicolumn{1}{c@{}}{$\bolds{\mathcal{M}_{\mathrm{union}}^{\dag,2}}$} \\
\hline
All correct & bias &   0.001  &   -0.207  &   0.498  & 0.003  &   -0.08 & -0.079  \\
  & $MC$ s.e.\tabnoteref[\ast]{vz} &    2.   614 &   8.   333  &  20.   214 &   2.   615 1 &    2.   615 5&    2.   615 3 \\
$Y$ wrong & bias &   -9.   87 & -10.   221 &   0.498  &   -0.147  & -0.502  &   -0.202  \\
  & $MC$ s.e. &   3.   322  & 10.   539 &   20.   214 & 4.   461  &   3.   177  & 3.   141  \\
$M$ wrong & bias & -0.033  &   -0.207  & -9.   497 &   0.001  &   0.046  &   0.046  \\
& $MC$ s.e. &    2.   613  & 8.   333  &  15.   376 & 2.   615  &    2.   614  &   2.   614  \\
$E$ wrong & bias &   -0.001  & 0.132 & 210.   450 & 0.066  & -0.089  & -0.087  \\
  & $MC$ s.e. & 2.614 & 4.373 & 2336.92 & 4.891 & 2.619 & 2.615 \\
$Y,E$ wrong & bias & -9.869 & -13.535 & 210.454 & -33.090 & -1.4609 & -2.487\\
  & $MC$ s.e. & 3.322 & 5.256 & 2336.92 & 375.334 & 5.187 & 4.245 \\
$Y,M$ wrong & bias & -9.355 & -10.220 & -9.496 & -4.346 & -3.579 &-3.579 \\
& $MC$ s.e. & 3.224 & 10.539 & 15.376 & 3.912 & 3.480 & 3.441 \\
$E,M$ wrong & bias & -0.032 & 0.132 & 205.060 & 0.088 & -0.001 & -3.77\times10^{-5} \\
  & $MC$ s.e. & 2.614 & 4.373 & 2289.788 & 4.763 & 2.623 & 2.618 \\
$Y,E,M$ wrong & bias & -9.355 & -13.535 & 205.060 & -37.757 & -4.223 & -5.253\\
  & $MC$ s.e. & 3.224 & 5.356 & 2289.78 & 379.122 & 5.835 & 4.828 \\
\hline
\end{tabular*}
\tabnotetext[]{tz}{$\mathcal{M}_{\mathrm{ym}}$: $\widehat{\theta}_{0}^{\,\mathrm{ym}}$; $\mathcal{M}_{\mathrm{ye}}$: $%
\widehat{\theta}_{0}^{\,\mathrm{ye}}$; $\mathcal{M}_{\mathrm{em}}$: $\widehat{\theta
}_{0}^{\,\mathrm{em}}$; $
\mathcal{M}_{\mathrm{union}}$: $\widehat{\theta}_{0}^{\,\mathrm{triply}}$; $\mathcal{M}%
_{\mathrm{union}}^{\dag,1}$: $\widehat{\theta}_{0}^{\,\mathrm{triply},\dag,1}$; $\mathcal{M}%
_{\mathrm{union}}^{\dag,2}$: $\widehat{\theta}_{0}^{\,\mathrm{triply},\dag,2}$.}
\tabnotetext[\ast]{vz}{Monte Carlo standard error.}
\end{table}

\begin{table}[t!]
\tabcolsep=0pt
\caption{Simulation results $n=1000$}\label{tab5}
\begin{tabular*}{\textwidth}{@{\extracolsep{\fill}}lld{2.8}d{3.3}d{2.6}d{5.6}d{2.4}d{2.4}@{}}
\hline
& & \multicolumn{1}{c}{$\bolds{\mathcal{M}_{\mathrm{ym}}}$} & \multicolumn{1}{c}{$\bolds{\mathcal{M}_{\mathrm{ye}}}$} &
\multicolumn{1}{c}{$\bolds{\mathcal{M}_{\mathrm{em}}}$} &
\multicolumn{1}{c}{$\bolds{\mathcal{M}_{\mathrm{union}}}$} & \multicolumn{1}{c}{$\bolds{\mathcal{M}_{\mathrm{union}}^{\dag,1}}$} &
\multicolumn{1}{c@{}}{$\bolds{\mathcal{M}_{\mathrm{union}}^{\dag,2}}$} \\
\hline
All correct & bias & 0.0324 & 0.004 & -0.106 & 0.034 & -0.047 & -0.047\\
$ $ & $MC$ s.e.\tabnoteref[\ast]{vz} & 1.136 & 3.06 & 6.490 & 1.136 & 1.137 &1.137 \\
$Y$ wrong & bias & -10.256 & -10.305 & -0.106 & 0.063 & -0.147 & -0.148\\
$ $ & $MC$ s.e. & 1.675 & 4.005 & 6.490 & 1.769 & 1.419 & 1.407 \\
$M$ wrong & bias & \multicolumn{1}{c}{$-5\times10^{-4}$\phantom{00.}} & 0.004 & -9.706 & 0.033 & 0.076 &0.076 \\
& $MC$ s.e. & 1.136 & 3.060 & 5.395 & 1.137 & 1.137 & 1.135 \\
$E$ wrong & bias & 0.032 & 0.135 & 2.4\times10^{6} & 1908.76 &-0.038 &-0.030 \\
$ $ & $MC$ s.e. & 1.136 & 1.794 & 4.3\times10^{7} & 53911.63 & 1.400 &1.242 \\
$Y,E$ wrong & bias & -10.256 & -14.011 & 2.4\times10^{6} &-1.1\times10^{6} & 6.201 & 1.024 \\
$ $ & $MC$ s.e. & 1.675 & 2.386 & 4.3\times10^{7} & 2.1\times10^{7} &9.406 & 5.097 \\
$Y,M$ wrong & bias & -9.705 & -10.305 & -9.706 & -4.216 & -3.555 &-3.557 \\
& $MC$ s.e. & 1.626 & 4.004 & 5.395 & 1.667 & 1.527 & 1.510 \\
$E,M$ wrong & bias & 5.7\times10^{-4} & 0.135 & 2.5\times10^{6} &2034.83 & 0.0539 & 0.0599 \\
$ $ & $MC$ s.e. & 1.136 & 1.794 & 4.6\times10^{7} & 56090.10 & 1.429 &1.272 \\
$Y,E,M$ wrong & bias & -9.075 & -14.011 & 2.5\times10^{6} &-1.2\times10^{6} & 4.659 & -0.755 \\
$ $ & $MC$ s.e. & 1.626 & 2.386 & 4.6\times10^{7} & 2.2\times10^{7} &10.121 & 5.910 \\
\hline
\end{tabular*}
\tabnotetext[]{tz}{$\mathcal{M}_{\mathrm{ym}}$: $\widehat{\theta}_{0}^{\,\mathrm{ym}}$; $\mathcal{M}_{\mathrm{ye}}$: $%
\widehat{\theta}_{0}^{\,\mathrm{ye}}$; $\mathcal{M}_{\mathrm{em}}$: $\widehat{\theta
}_{0}^{\,\mathrm{em}}$; $
\mathcal{M}_{\mathrm{union}}$: $\widehat{\theta}_{0}^{\,\mathrm{triply}}$; $\mathcal{M}%
_{\mathrm{union}}^{\dag,1}$: $\widehat{\theta}_{0}^{\,\mathrm{triply},\dag,1}$; $\mathcal{M}%
_{\mathrm{union}}^{\dag,2}$: $\widehat{\theta}_{0}^{\,\mathrm{triply},\dag,2}$.}
\tabnotetext[\ast]{vz}{Monte Carlo standard error.}
\end{table}

Tables~\ref{tab4} and~\ref{tab5} summarize simulation results for $\widehat{\theta
}_{0}^{\,\mathrm{ym}}$, $
\widehat{\theta}_{0}^{\,\mathrm{ye}}$, $\widehat{\theta}_{0}^{\,\mathrm{em}}$, $\widehat{\theta}%
_{0}^{\,\mathrm{triply}}$, $\widehat{\theta}_{0}^{\,\mathrm{triply},\dag,1}$ and $\widehat{%
\theta}_{0}^{\,\mathrm{triply},\dag,2}$. When all three working models are correct,
all estimators perform well in terms of bias, but there are clear
differences between the estimators in terms of efficiency. In fact, $%
\widehat{\theta}_{0}^{\,\mathrm{ym}}$, $\widehat{\theta}_{0}^{\,\mathrm{triply}}$, $\widehat
{\theta
}_{0}^{\,\mathrm{triply},\dag,1}$ and $\widehat{\theta}_{0}^{\,\mathrm{triply},\dag,2}$ have
comparable efficiency for $n=200,1000$, but $\widehat{\theta}_{0}^{\,\mathrm{ye}}$,
$\widehat{\theta}_{0}^{\,\mathrm{em}}$ is far more variable. Moreover, under
mis-specification of a single model, $\widehat{\theta}_{0}^{\,\mathrm{triply}}$, $
\widehat{\theta}_{0}^{\,\mathrm{triply},\dag,1}$ and $\widehat{\theta}%
_{0}^{\,\mathrm{triply},\dag,2}$ remain nearly unbiased, and for the most part
substantially more efficient than the corresponding consistent
estimator in $%
\{\widehat{\theta}_{0}^{\,\mathrm{ym}}$, $\widehat{\theta}_{0}^{\,\mathrm{ye}}$, $\widehat{\theta
}%
_{0}^{\,\mathrm{em}}\}$. When at least two models are mis-specified, the multiply
robust estimators $\widehat{\theta}_{0}^{\,\mathrm{triply}}$, $\widehat{\theta}%
_{0}^{\,\mathrm{triply},\dag,1}$ and $\widehat{\theta}_{0}^{\,\mathrm{triply},\dag,2}$
generally outperform the other estimators, although $\widehat{\theta}%
_{0}^{\,\mathrm{triply}}$ occasionally succumbs to the unstable weights resulting in
disastrous mean squared error; see Table~\ref{tab5} when Model~M and Model~E are both
incorrect. In contrast, $\widehat{\theta}_{0}^{\,\mathrm{triply},\dag,2}$ generally
improves on $\widehat{\theta}_{0}^{\,\mathrm{triply},\dag,1}$ which generally
outperforms $\widehat{\theta}_{0}^{\,\mathrm{triply}}$ and for the most part
$\widehat{%
\theta}_{0}^{\,\mathrm{triply},\dag,1}$ and $\widehat{\theta}_{0}^{\,\mathrm{triply},\dag,2}$
appear to eliminate any possible deleterious impact of highly variable
weights.

\section{A comparison to some existing estimators}\label{sec7}

In this section, we briefly compare the proposed approach to some existing
estimators in the literature. Perhaps the most common approach for
estimating direct and indirect effects when $Y$ is continuous uses a system
of linear structural equations; whereby, a linear structural equation for
the outcome, given the exposure, the mediator and the confounders, is combined
with a linear structural equation for the mediator, given the exposure and
confounders, to produce an estimator of natural direct and indirect effects.
The classical approach of \citet{Baron} is a particular instance of
this approach. In recent work, mainly motivated by Pearl's mediation
functional, several authors [Imai, Keele and Tingley (\citeyear{Imai2}), Imai,
Keele and Yamamoto (\citeyear{Imai1}), Pearl (\citeyear{Pearl2}),
\citet{vanderweele2}, \citet{Vanderweele}] have
demonstrated how the simple
linear structural equation approach generalizes to accommodate both, the
presence of an interaction between exposure and mediator variables, and a
nonlinear link function, either in the regression model for the
outcome, or in
the regression model for the mediator, or both. In fact, when the
effect of
confounders is also modeled in such structural equations, inferences based
on the latter can be viewed as special instances of inferences obtained
under a particular specification of model $\mathcal{M}_{a}$ for the outcome
and the mediator densities. And thus, as previously shown in the
simulations, an estimator obtained under a system of structural equations
will generally fail to produce a consistent estimator of natural direct and
indirect effects when model $\mathcal{M}_{a}$ is incorrect, whereas, by using
the proposed multiply robust estimator, valid inferences can be recovered
under the union model $\mathcal{M}_{b}\cup\mathcal{M}_{c}$, even if $%
\mathcal{M}_{a}$ fails.

A notable improvement on the system of structural equations approach is the
double robust estimator of a natural direct effect due to \citet{vanderlaan2}. Their estimator solves the estimating equation
constructed using an empirical version of $S_{\mathrm{NDE},\mathrm{singleton}}^{\mathrm{eff},%
\mathcal{M}_{a}\cup\mathcal{M}_{c}} ( \theta_{0},\delta_{0} ) $
given in the online Appendix. They show their estimator remains CAN in
the larger submodel $\mathcal{M}_{a}\cup\mathcal{M}_{c}$ and therefore,
they can recover valid inferences even when the outcome model is incorrect,
provided both the exposure and mediator models are correct. Unfortunately,
the van der Laan estimator is still not entirely satisfactory because unlike
the proposed multiply robust estimator, it requires that the model for the
mediator density is correct. Nonetheless, if the\vadjust{\goodbreak} mediator model is correct,
the authors establish that their estimator achieves the efficiency
bound for
model $\mathcal{M}_{a}\cup\mathcal{M}_{c}$ at the intersection
submodel $\mathcal{M}_{a}\cap\mathcal{M}_{c}$ where all models are correct; and thus
it is locally semiparametric efficient in $\mathcal{M}_{a}\cup\mathcal
{M}_{c}$. Interestingly, as we report in the online supplement, the
semiparametric efficiency bounds for models $\mathcal{M}_{a}\cup
\mathcal{M}_{c}$ and $\mathcal{M}_{a}\cup\mathcal{M}_{b}\cup\mathcal{M}_{c}$ are
distinct, because the density of the mediator variable is not ancillary for
inferences about the M-functional. Thus, any restriction placed on the
mediator's conditional density can, when correct, produce improvements in
efficiency. This is in stark contrast with the role played by the density
of the exposure variable, which as in the estimation of the marginal causal
effect, remains ancillary for inferences about the M-functional and
thus the
efficiency bound for the latter is unaltered by any additional information
on the former [\citet{Robins6}]. In the online Appendix, we provide a
general functional map that relates the efficient influence function
for the
larger model $\mathcal{M}_{a}\cup\mathcal{M}_{b}\cup\mathcal{M}_{c}$ to
the efficient influence for the smaller model $\mathcal{M}_{a}\cup
\mathcal{M}_{c}$ where the model for the mediator is either parametric or
semiparametric. Our map is instructive because it makes explicit using
simple geometric arguments, the information that is gained from increasing
restrictions on the law of the mediator. In the online Appendix, we
illustrate the map by recovering the efficient influence function of
van der
Laan and Petersen in the case of a singleton model (i.e., a~known conditional
density) for
the mediator and in the case of a parametric model for the mediator.

\section{A semiparametric sensitivity analysis}\label{sec8}

We describe a semiparametric sensitivity analysis framework to assess the
extent to which a violation of the ignorability assumption for the mediator
might alter inferences about natural direct and indirect effects. Although
only results for the natural direct effect are given here, the
extension for
the indirect effect is easily deduced from the presentation. Let
\[
t ( e,m,x ) =\mathbb{E} [ Y_{1,m}|E=e,M=m,X=x ] -\mathbb{E}%
[ Y_{1,m}|E=e,M\neq m,X=x ] ,
\]
then
\[
Y_{e^{\prime},m}\not\perp\!\!\!\perp M|E=e,X,
\]
that is, a violation of the ignorability assumption for the mediator variable,
generally implies that
\[
t ( e,m,x ) \neq0\qquad\mbox{for some }(e,m,x).
\]
Thus, we proceed as in \citet{Robins22}, and propose
to recover inferences by assuming the selection bias function $t (
e,m,x ) $ is known, which encodes the magnitude and direction of the
unmeasured confounding for the mediator. In the following, the support of $M$, $\mathcal{S}$ is
assumed to be finite. To motivate the proposed approach, suppose for the
moment that\vadjust{\goodbreak} $f_{M|E,X} ( M|E,X ) $ is known; then under the
assumption that the exposure is ignorable given $X$, we show in the \hyperref[app]{Appendix}
that
\begin{eqnarray*}
&&\mathbb{E} [ Y_{1,m}|M_{0}=m,X=x ]\\
&&\qquad =\mathbb{E} [
Y_{1,m}|E=0,M=m,X=x ] \\
&&\qquad=\mathbb{E} [ Y|E=1,M=m,X=x ] -t ( 1,m,x )
\bigl( 1-f_{M|E,X} ( m|E=1,X=x ) \bigr) \\
&&\quad\qquad{}+t ( 0,m,x ) \bigl( 1-f_{M|E,X} ( m|E=0,X=x ) \bigr) ,
\end{eqnarray*}
and therefore the M-functional is identified by%
%
\begin{eqnarray}\label{rep1}
&&\sum_{m\in\mathcal{S}}\mathbb{E} \bigl\{ \mathbb{E} [ Y|E=1,M=m,X%
] -\mbox{ }t ( 1,m,X ) \bigl( 1-f_{M|E,X} (
m|E=1,X ) \bigr) \nonumber\\
&&\hspace*{137pt}{} +t ( 0,m,X ) \bigl( 1-f_{M|E,X} ( m|E=0,X )
\bigr) \bigr\}\\
&&\qquad{}\times f_{M|E,X} ( m|E=0,X ) , \nonumber
\end{eqnarray}
which is equivalently represented as
%
\begin{eqnarray} \label{rep2}
&&\mathbb{E} \biggl[ \frac{I \{ E=1 \} f_{M|E,X} ( M|E=0,X )
}{f_{E|X}(1|X)f_{M|E,X} ( M|E=1,X ) } \nonumber\\
&&\hspace*{10pt}{}\times \bigl\{ Y-t ( 1,M,X ) \bigl(
1-f_{M|E,X} ( m|E=1,X ) \bigr)\\
&&\hspace*{33pt}{} +t ( 0,M,X ) \bigl(
1-f_{M|E,X} ( M|E=0,X ) \bigr) \bigr\} \biggr] . \nonumber
\end{eqnarray}
Below, these two equivalent representations,  (\ref{rep1})  and
(\ref{rep2}), are carefully combined to obtain a double robust
estimator of the M-functional, assuming $t ( \cdot,\cdot,\cdot ) $
is known. A sensitivity analysis is then obtained by repeating this process
and reporting inferences for each choice of $t ( \cdot,\cdot,\cdot
) $ in a finite set of user-specified functions $\mathcal{T}=\{$ $%
t_{\lambda} ( \cdot,\cdot,\cdot ) \dvtx\lambda\}$ indexed by a
finite dimensional parameter $\lambda$ with $t_{0} ( \cdot,\cdot
,\cdot ) \in\mathcal{T}$ corresponding to the unmeasured
confounding assumption, that is, $t_{0} ( \cdot,\cdot,\cdot ) \equiv
0$. Throughout, the model $f_{M|E,X}^{\,\mathrm{par}}(\cdot|E,X;\beta_{m})$ for the
probability mass function of $M$ is assumed to be correct. Thus, to
implement the sensitivity analysis, we develop a semiparametric
estimator of
the natural direct effect in the union model $\mathcal{M}_{a}\cup
\mathcal{M}%
_{c} $, assuming $t ( \cdot,\cdot,\cdot ) $ =$t_{\lambda^{\ast
}} ( \cdot,\cdot,\cdot ) $ for a fixed $\lambda^{\ast}$. The
proposed doubly robust estimator of the natural direct effect is then given
by $\widehat{\theta}_{0}^{\,\mathrm{doubly}} ( \lambda^{\ast} ) -\widehat{%
\delta}_{0}^{\,\mathrm{doubly}}$ where $\widehat{\delta}_{0}^{\,\mathrm{doubly}}$ is as
previously described, and
\begin{eqnarray*}
\widehat{\theta}_{0}^{\,\mathrm{doubly}} ( \lambda^{\ast} ) &=&\mathbb{P}%
_{n} \biggl[ \frac{I \{ E=1 \} \widehat{f}_{M|E,X}^{\,\mathrm{par}} (
M|E=0,X ) }{\widehat{f}_{E|X}^{\,\mathrm{par}}(1|X)\widehat{f}_{M|E,X}^{\,\mathrm{par}} (
M|E=1,X ) }\\
&&\hspace*{14pt}{}\times \{ Y-\widehat{\mathbb{E}}^{\mathrm{par}} ( Y|X,M,E=1 )
 \}
 +\widetilde{\eta}^{\,\mathrm{par}} ( 1,0,X;\lambda^{\ast} ) \biggr],
\end{eqnarray*}
with%
\begin{eqnarray*}
&&\hspace*{-5pt}\widetilde{\eta}^{\,\mathrm{par}} ( 1,0,X;\lambda^{\ast} )\\[-2pt]
 &&\hspace*{-5pt}\quad=\sum_{m\in
\mathcal{S}} \bigl\{ \widehat{\mathbb{E}}^{\mathrm{par}} ( Y|X,M=m,E=1 )
+t_{\lambda^{\ast}} ( 0,m,X ) \bigl( 1-\widehat{f}%
_{M|E,X}^{\,\mathrm{par}} ( m|E=0,X ) \bigr) \\[-2pt]
&&\hspace*{-5pt}\hspace*{165pt}{} -t_{\lambda^{\ast}} ( 1,m,X ) \bigl( 1-\widehat{f}%
_{M|E,X}^{\,\mathrm{par}} ( m|E=1,X ) \bigr) \bigr\}\\[-2pt]
&&\hspace*{-5pt}\hspace*{17pt}\qquad{}\times \widehat{f}%
_{M|E,X}^{\,\mathrm{par}} ( m|E=0,X ) .
\end{eqnarray*}

  Our sensitivity analysis then entails reporting the set $ \{
\widehat{\theta}_{0}^{\,\mathrm{doubly}} ( \lambda ) -\break\widehat{\delta}%
_{0}^{\,\mathrm{doubly}}\dvtx\lambda \} $ (and the associated confidence intervals),
which summarizes how sensitive inferences are to a deviation from the
ignorability assumption $\lambda=0$. A theoretical justification for the
approach is given by the following formal result, which is proved in the
supplemental Appendix.

\renewcommand{\thetheorem}{\arabic{theorem}}
\setcounter{theorem}{3}
\begin{theorem}\label{teo4}
 Suppose $t ( \cdot,\cdot,\cdot ) =
t_{\lambda^{\ast}} ( \cdot,\cdot,\cdot )$;
then under the consistency, positivity assumptions and the ignorability
assumption for the exposure,\break $\widehat{\theta}_{0}^{\,\mathrm{doubly}} ( \lambda
^{\ast} ) -\widehat{\delta}_{0}^{\,\mathrm{doubly}}$ is a CAN estimator
of the natural direct effect in $\mathcal{M}_{a}\cup\mathcal{M}_{c}$.
\end{theorem}

The influence function of $\widehat{\theta}_{0}^{\,\mathrm{doubly}} ( \lambda
^{\ast} ) $ is provided in the \hyperref[app]{Appendix}, and can be used to construct
a corresponding confidence interval.

It is important to note that the sensitivity analysis technique presented
here differs in crucial ways from previous techniques developed by
\citet{Hafeman2}, \citet{Vandeweele3} and \citet{Imai1}. First, the methodology of
\citet{Vandeweele3} postulates the existence of an unmeasured confounder $U$
(possibly vector valued) which, when included in~$X$, recovers the sequential
ignorability assumption. The sensitivity analysis then requires
specification of a sensitivity parameter encoding the effect of the
unmeasured confounder on the outcome within levels of $( E,X,M )$, and another parameter for the effect of the exposure on the density
of the
unmeasured confounder given $(X,M)$. This is a daunting task which renders
the approach generally impractical, except perhaps in the simple setting
where it is reasonable to postulate a single binary confounder is
unobserved, and one is willing to make further simplifying assumptions about
the required sensitivity parameters [\citet{Vandeweele3}]. In
comparison, the
proposed approach circumvents this difficulty by concisely encoding a
violation of the ignorability assumption for the mediator through the
selection bias function $t_{\lambda} ( e,m,x ) $. Thus the
approach makes no reference and thus is agnostic about the existence,
dimension and nature of unmeasured confounders $U$. Furthermore, in our
proposal, the ignorability violation can arise due to an unmeasured
confounder of the mediator-outcome relationship that is also an effect of
the exposure variable, a setting not handled by the technique of
\citet{Vandeweele3}. The method of \citet{Hafeman2} which is\vadjust{\goodbreak}
restricted to binary data,
shares some of the limitations given above. Finally, in contrast with our
proposed double robust approach, a coherent implementation of the
sensitivity analysis techniques of \citet{Imai1}, \citet{Imai2} and
\citet{Vandeweele3} rely on correct specification of all
posited models. We refer
the reader to \citet{Vandeweele3} for further discussion of
\citet{Hafeman2}
and \mbox{\citet{Imai1}}.

\section{ Discussion}\label{sec9}

The main contribution of the current paper is a theoretically rigorous yet
practically relevant semiparametric framework for making inferences about
natural direct and indirect causal effects in the presence of a large number
of confounding factors. Semiparametric efficiency bounds are given for the
nonparametric model, and multiply robust locally efficient estimators are
developed that can be used when nonparametric estimation is not possible.

Although the paper focuses on a binary exposure, we note that the extension
to a polytomous exposure is trivial. In future work, we shall extend our
results for marginal effects by considering conditional natural direct and
indirect effects, given a subset of pre-exposure variables [Tchetgen Tchetgen
and Shpitser (\citeyear{TchShpi11})]. These
models are
particularly important in making inferences about so-called moderated
mediation effects, a topic of growing interest, particularly in the
field of
psychology [\citet{Preacher}]. In related work, we have
recently extended our results to a survival analysis setting [\citet{Tchetgen}].

A major limitation of the current paper is that it assumes that the mediator
is measured without error, an assumption that may be unrealistic in
practice and, if incorrect, may result in biased inferences about mediated
effects. We note that much of the recent literature on causal mediation
analysis makes a similar assumption. In future work, it will be
important to
build on the results derived in the current paper to appropriately account
for a mis-measured mediator [Tchetgen Tchetgen
and Lin (\citeyear{Tch12})].

\begin{appendix}\label{app}
\section*{Appendix}

\begin{pf*}{Proof of Theorem~\ref{teo1}}
Let $F_{O;t} =F_{Y|M,X,E;t}F_{M|E,X;t}F_{E|X;t}F_{X;t}$ denote a
one-dimensional regular parametric submodel of $\mathcal{M}_{\mathrm{nonpar}}$,
with $F_{O,0}=F_{O}$, and let
\begin{eqnarray*}
&&\theta_{t}=\theta_{0} ( F_{O;t} ) =\mathop{\int\hspace*{-4pt}\int}_{\mathcal{S\times X}}
\mathbb{E}_{t} ( Y|E=1,M=m,X=x )\\
&&\hspace*{91pt}{}\times f_{M|E,X;t} ( m|E=0,X=x ) f_{X;t}(x)\,d\mu(m,x).
\end{eqnarray*}
The efficient influence function $S_{\theta_{0}}^{\mathrm{eff},\mathrm{nonpar}} ( \theta
_{0} ) $ is the unique random variable to satisfy the following
equation:
\[
\nabla_{t=0}\theta_{t}=\mathbb{E} \{ S_{\theta
_{0}}^{\mathrm{eff},\mathrm{nonpar}} ( \theta_{0} ) U \}\vadjust{\goodbreak}
\]
for $U$ the score of $F_{O;t}$ at $t=0$, and $\nabla_{t=0}$ denoting
differentiation w.r.t. $t$ at $t=0$. We observe that
\begin{eqnarray*}
\frac{\partial\theta_{t}}{\partial t}\bigg|_{t=0} &=&\mathop{\int\hspace*{-4pt}\int}_{\mathcal{S\times X}}
\nabla_{t=0}\mathbb{E}_{t} (
Y|E=1,M=m,X=x )\\
&&\hspace*{20pt}{}\times f_{M|E,X} ( m|E=0,X=x ) f_{X}(x)\,d\mu(m,x) \\
&&{}+\mathop{\int\hspace*{-4pt}\int}_{\mathcal{S\times X}}\mathbb{E}%
( Y|E=1,M=m,X=x )\\
&&\hspace*{29pt}{}\times \nabla_{t=0}f_{M|E,X;t} ( m|E=0,X=x )
f_{X}(x)\,d\mu(m,x) \\
&&{}+\mathop{\int\hspace*{-4pt}\int}_{\mathcal{S\times X}}\mathbb{E}%
( Y|E=1,M=m,X=x )\\
&&\hspace*{29pt}{}\times f_{M|E,X} ( m|E=0,X=x ) \nabla
_{t=0}f_{X;t}(x)\,d\mu(m,x).
\end{eqnarray*}
Considering the first term, it is straightforward to verify that
\begin{eqnarray*}
&&\mathop{\int\hspace*{-4pt}\int}_{\mathcal{S\times X}}\nabla_{t=0}%
\mathbb{E}_{t} ( Y|E=1,M=m,X=x ) f_{M|E,X} ( m|E=0,X=x )
f_{X}(x)\,d\mu(m,x) \\
&&\qquad=\mathbb{E} \biggl[ U\frac{I(E=1)}{f_{E|X} ( E|X ) } \{ Y-%
\mathbb{E} ( Y|E,M=m,X=x ) \} \frac{f_{M|E,X} (
M|E=0,X ) }{f_{M|E,X} ( M|E=1,X ) } \biggr].
\end{eqnarray*}
Similarly, one can easily verify that
\begin{eqnarray*}
&&\hspace*{-4pt}\mathop{\int\hspace*{-4pt}\int}_{\mathcal{S\times X}}\mathbb{E}%
( Y|E=1,M=m,X=x ) \nabla_{t=0}f_{M|E,X;t} ( m|E=0,X=x )
f_{X}(x)\,d\mu(m,x) \\
&&\qquad=\mathbb{E} \biggl[ U\frac{I(E=0)}{f_{E|X} ( E|X ) } \{ \mathbb{%
E} ( Y|E=1,M=m,X=x ) -\eta ( 1,0,X ) \} \biggr],
\end{eqnarray*}
and finally, one can also verify that
\begin{eqnarray*}
&&\hspace*{-4pt}\mathop{\int\hspace*{-4pt}\int}_{\mathcal{S\times X}}\mathbb{E}%
( Y|E=1,M=m,X=x ) f_{M|E,X} ( m|E=0,X=x ) \nabla
_{t=0}f_{X;t}(x)\,d\mu(m,x) \\
&&\qquad=\mathbb{E} [ U \{ \eta ( 1,0,X ) -\theta_{0} \} %
].
\end{eqnarray*}
Thus we obtain
\[
\nabla_{t=0}\theta_{t}=\mathbb{E} \{ S_{\theta
_{0}}^{\mathrm{eff},\mathrm{nonpar}} ( \theta_{0} ) U \}.
\]
Given $S_{\delta e}^{\mathrm{eff},\mathrm{nonpar}} ( \delta_{e} ) $, the results for
the direct and indirect effect follow from the fact that the influence
function of a difference of two functionals equals the difference of the
respective influence functions. Because the model is nonparametric, there
is a unique influence function for each functional, and it is efficient in
the model, leading to the efficiency bound results.\vadjust{\goodbreak}
\end{pf*}

\begin{pf*}{Proof of Theorem~\ref{teo2}}
We begin by showing that
\begin{eqnarray}\label{appzero}
&& \mathbb{E}\{S_{\theta_{0}}^{\mathrm{eff},\mathrm{nonpar}} ( \theta_{0};\beta
_{m}^{\ast},\beta_{e}^{\ast},\beta_{y}^{\ast} ) \}
 \nonumber\\[-8pt]\\[-8pt]
&&\qquad =0 \nonumber
\end{eqnarray}
under model $\mathcal{M}_{\mathrm{union}}$. First note that $ ( \beta
_{y}^{\ast},\beta_{m}^{\ast} ) = ( \beta_{y},\beta_{m} ) $
under model $\mathcal{M}_{a}$. Equality (\ref{appzero}) now follows
because $\mathbb{E}^{\mathrm{par}} ( Y|X,M,E=1;\beta_{y} ) =\mathbb{E} (
Y|X,M,E=1 ) $ and $\eta ( 1,0,X;\beta_{y},\beta_{m} ) =\mathbb{E} [ \{ \mathbb{E}^{\mathrm{par}} ( Y|X,M,E=1;\beta_{y} )
\} |E=0,X ] =\eta ( 1,0,X )$:
\begin{eqnarray*}
& &\mathbb{E}\{S_{\theta_{0}}^{\mathrm{eff},\mathrm{nonpar}} ( \theta_{0};\beta
_{m},\beta_{e}^{\ast},\beta_{y} ) \} \\
&&\qquad =\mathbb{E} \biggl[ \frac{I \{ E=1 \} f_{M|E,X}^{\,\mathrm{par}} (
M|E=0,X;\beta_{m} ) }{f_{E|X}^{\,\mathrm{par}}(1|X;\beta_{e}^{\ast
})f_{M|E,X}^{\,\mathrm{par}} ( M|E=1,X;\beta_{m} ) }\\
&&\hspace*{42pt}{}\times{\overbrace{\mathbb{E} \{ Y-\mathbb{E}^{\mathrm{par}} ( Y|X,M,E=1;\beta_{y} )
|E=1,M,X \}}^ {\mathrm{=0}}} \biggr] \\
&&\quad\qquad{} +\mathbb{E} \biggl[
\frac{I(E=0)}{f_{E|X}^{\,\mathrm{par}}(1|X;\beta_{e}^{\ast})}\\
&&\hspace*{55pt}{}\times{\overbrace{\mathbb{E} [ \{ \mathbb{E}^{\mathrm{par}} (
Y|X,M,E=1;\beta_{y} ) -\eta ( 1,0,X;\beta_{y},\beta_{m} )
\} |E=0,X ] }^{\mathrm{=0}}} \biggr] \\
&&\quad\qquad{} +\mathbb{E} [ \eta ( 1,0,X;\beta_{y},\beta_{m} ) ]
-\theta_{0} \\
&&\qquad =0.
\end{eqnarray*}
Second, $ ( \beta_{y}^{\ast},\beta_{e}^{\ast} ) = (
\beta_{y},\beta_{e} ) $ under model $\mathcal{M}_{b}$. Equality (\ref%
{appzero}) now follows because $\mathbb{E}^{\mathrm{par}} ( Y|X,M,E=1;\beta
_{y} ) =\mathbb{E} ( Y|X,M,E=1 ) $ and $f_{E|X}^{\,\mathrm{par}}(1|X;%
\beta_{e})=\break f_{E|X}(1|X)$:
\begin{eqnarray*}
&&\mathbb{E}\{S_{\theta_{0}}^{\mathrm{eff},\mathrm{nonpar}} ( \theta_{0};\beta
_{m}^{\ast},\beta_{e},\beta_{y} ) \} \\
&&\qquad=\mathbb{E} \biggl[ \frac{I \{ E=1 \} f_{M|E,X}^{\,\mathrm{par}} (
M|E=0,X;\beta_{m}^{\ast} ) }{f_{E|X}^{\,\mathrm{par}}(1|X;\beta
_{e})f_{M|E,X}^{\,\mathrm{par}} ( M|E=1,X;\beta_{m}^{\ast} ) }\\
&&\hspace*{43pt}{}\times{%
\overbrace{\mathbb{E} \{ Y-\mathbb{E}^{\mathrm{par}} ( Y|X,M,E=1;\beta
_{y} ) |E=1,M,X \} }^{\mathrm{=0}}} \biggr] \\
&&\quad\qquad{}+\mathbb{E} \biggl[ \frac{I(E=0)}{f_{E|X}^{\,\mathrm{par}}(1|X;\beta_{e})}\\
&&\hspace*{56pt}{}\times\mathbb{E}%
[ \{ \mathbb{E}^{\mathrm{par}} ( Y|X,M,E=1;\beta_{y} ) -\eta
( 1,0,X;\beta_{y},\beta_{m}^{\ast} ) \} |E=0,X ] %
\biggr] \\
&&\quad\qquad{}+\mathbb{E} [ \eta ( 1,0,X;\beta_{y},\beta_{m}^{\ast} ) %
] -\theta_{0} \\
&&\qquad=\mathbb{E} \bigl[ \mathbb{E} [ \{ \mathbb{E}^{\mathrm{par}} (
Y|X,M,E=1;\beta_{y} ) \} |E=0,X ] \bigr] -\theta_{0}=0.
\end{eqnarray*}
Third, equality (\ref{appzero}) holds under model $\mathcal{M}_{c}$
because
\begin{eqnarray*}
&&\mathbb{E}\{S_{\theta_{0}}^{\mathrm{eff},\mathrm{nonpar}} ( \theta_{0};\beta_{m}^{\
},\beta_{e} ,\beta_{y}^{\ast} ) \} \\
&&\qquad=\mathbb{E} \biggl[ \frac{I \{ E=1 \} f_{M|E,X}^{\,\mathrm{par}} (
M|E=0,X;\beta_{m} ) }{f_{E|X}^{\,\mathrm{par}}(1|X;\beta
_{e})f_{M|E,X}^{\,\mathrm{par}} ( M|E=1,X;\beta_{m} ) }\\
&&\hspace*{43pt}{}\times\mathbb{E} \{ Y-%
\mathbb{E}^{\mathrm{par}} ( Y|X,M,E=1;\beta_{y}^{\ast} )  \} \biggr]
\\
&&\quad\qquad{}+\mathbb{E} \biggl[ \frac{I(E=0)}{f_{E|X}^{\,\mathrm{par}}(1|X;\beta_{e})}\\
&&\hspace*{57pt}{}\times\mathbb{E}%
[ \{ \mathbb{E}^{\mathrm{par}} ( Y|X,M,E=1;\beta_{y}^{\ast} )
-\eta ( 1,0,X;\beta_{y}^{\ast},\beta_{m} ) \} |E=0,X%
] \biggr] \\
&&\quad\qquad{}+\mathbb{E} [ \eta ( 1,0,X;\beta_{y}^{\ast},\beta_{m} ) %
] -\theta_{0} \\
&&\qquad=\mathbb{E} \bigl[ \mathbb{E} [ \{ \mathbb{E} (
Y|X,M,E=1 ) \} |E=0,X ] \bigr]\\
&&\quad\qquad{} -\mathbb{E} \bigl[ \mathbb{E}%
[ \mathbb{E}^{\mathrm{par}} ( Y|X,M,E=1;\beta_{y}^{\ast} ) |E=0,X%
] \bigr] \\
&&\quad\qquad{}+\mathbb{E} \bigl[ \mathbb{E} [ \mathbb{E}^{\mathrm{par}} ( Y|X,M,E=1;\beta
_{y}^{\ast} ) |E=0,X ] \bigr] -\mathbb{E} [ \eta (
1,0,X;\beta_{y}^{\ast},\beta_{m} ) ] \\
&&\quad\qquad{}+\mathbb{E} [ \eta ( 1,0,X;\beta_{y}^{\ast},\beta_{m} ) %
] -\theta_{0} \\
&&\qquad=\mathbb{E} \bigl[ \mathbb{E} [ \{ \mathbb{E} (
Y|X,M,E=1 ) \} |E=0,X ] \bigr] -\theta_{0}.
\end{eqnarray*}

Assuming that the regularity conditions of Theorem 1A in \citet
{robins11} hold for $S_{\theta_{0}}^{\mathrm{eff},\mathrm{nonpar}} ( \theta_{0};\beta
_{m} ,\beta_{e} ,\beta_{y}  ) ,S_{\beta} ( \beta
) $, the expression for $S_{\theta_{0}}^{\mathrm{union}} ( \theta
_{0},\beta^{\ast} ) $ follows by standard Taylor expansion arguments,
and it now follows that
%
\begin{equation}\label{al}
\sqrt{n}(\widehat{\theta}_{0}^{\,\mathrm{triply}}-\theta_{0})=\frac{1}{n^{1/2}}%
\sum_{i=1}^{n}S_{\theta_{0},i}^{\mathrm{union}} ( \theta_{0},\beta^{\ast
} ) +o_{p}(1) .
\end{equation}
The asymptotic distribution of $\sqrt{n}(\widehat{\theta}%
_{0}^{\,\mathrm{triply}}-\theta_{0})$ under model $\mathcal{M}_{\mathrm{union}}$ follows
from the previous equation by Slutsky's Theorem and the Central Limit
Theorem.

We note that $\widehat{\delta}_{e}^{\,\mathrm{doubly}}$ is CAN in the union model
$%
\mathcal{M}_{\mathrm{union}}$ since it is CAN in the larger model where either
the density for the exposure is correct, or the density of the mediator and
the outcome regression are both correct and thus $\eta ( e,e,X;\beta
_{y}^{\ast},\beta_{m}^{\ast} ) =\mathbb{E} ( Y|X,E=e ) $.
This gives the multiply robust result for direct and indirect effects. The
asymptotic distribution of direct and indirect effect estimates then follows
from similar arguments as above.

At the intersection submodel
\[
\frac{\partial\mathbb{E} \{ S_{\theta_{0}}^{\mathrm{eff},\mathrm{nonpar}} ( \theta
_{0},\beta ) \} }{\partial\beta^{T}}=0
\]
hence
\[
S_{\theta_{0}}^{\mathrm{union}} ( \theta_{0},\beta ) =S_{\theta
_{0}}^{\mathrm{eff},\mathrm{nonpar}} ( \theta_{0},\beta ) .
\]
The semiparametric efficiency claim then follows for $\widehat{\theta}%
_{0}^{\,\mathrm{triply}}$, and a similar argument gives the result for direct and
indirect effects.
\end{pf*}

\begin{pf*}{Proofs of Theorems 3 and~\ref{teo4}}
The proofs are given in the online
Appendix.
\end{pf*}

\end{appendix}

\section*{Acknowledgments}
The authors would like to acknowledge Andrea Rotnitzky who provided
invaluable comments that improved the presentation of the results given in
Section~\ref{sec7}. The authors also thank James Robins and Tyler
VanderWeele for useful comments that significantly improved the presentation
of this article.

\begin{supplement}
\stitle{Supplemental Appendix to Semiparametric theory for causal mediation analysis}
\slink[doi]{10.1214/12-AOS990SUPP} 
\sdatatype{.pdf}
\sfilename{aos990\_supp.pdf}
\sdescription{The supplementary material gives
the semiparametric efficiency theory for estimation of natural direct effects with
a known model for the mediator density. The Appendix also gives the proof of
Theorem 3 (stated in the Supplementary Appendix) and of Theorem~\ref{teo4}.}
\end{supplement}

%

\printaddresses

\end{document}